\UseRawInputEncoding 

\documentclass{article}

\newtheorem{theorem}{Theorem}[section]
\newtheorem{definition}[theorem]{Definition}
\newtheorem{proposition}[theorem]{Proposition}
\newtheorem{lemma}[theorem]{Lemma}

\usepackage{proof}
\usepackage{latexsym}
\usepackage{amsfonts}
\usepackage{amssymb}
\usepackage{cite}

\newcommand{\Rarw }{\Rightarrow}
\newcommand{\Lrarw}{\Leftrightarrow}
\newcommand{\lrarw}{\leftrightarrow}
\newcommand{\benu}{\begin{enumerate}}
\newcommand{\eenu}{\end{enumerate}}
\newcommand{\restrict}{\!\upharpoonright\!}
\newcommand{\spand}{\,\&\,}

\newcommand{\supp}{{\rm supp}}

\begin{document}

\title{Wellfoundedness proof with the maximal distinguished set
}

\author{Toshiyasu Arai
\\
Graduate School of Mathematical Sciences
\\
University of Tokyo
\\
3-8-1 Komaba, Meguro-ku,
Tokyo 153-8914, JAPAN
\\
tosarai@ms.u-tokyo.ac.jp
}
\date{}

\maketitle

\begin{abstract}
In \cite{singlestable} it is shown that an ordinal
$\sup_{N<\omega}\psi_{\Omega_{1}}(\varepsilon_{\Omega_{\mathbb{S}+N}+1})$ is an upper bound for the proof-theoretic ordinal of a set theory ${\sf KP}\ell^{r}+(M\prec_{\Sigma_{1}}V)$.
In this paper we show that 
a second order arithmetic $\Sigma^{1-}_{2}\mbox{{\rm -CA}}+\Pi^{1}_{1}\mbox{{\rm -CA}}_{0}$
 proves the wellfoundedness up to $\psi_{\Omega_{1}}(\varepsilon_{\Omega_{\mathbb{S}+N+1}})$ for each $N$.
It is easy to interpret $\Sigma^{1-}_{2}\mbox{{\rm -CA}}+\Pi^{1}_{1}\mbox{{\rm -CA}}_{0}$
in ${\sf KP}\ell^{r}+(M\prec_{\Sigma_{1}}V)$.
\end{abstract}

\section{Introduction}\label{sect:introduction}
In \cite{singlestable} the following theorem is shown, where
${\sf KP}\ell^{r}+(M\prec_{\Sigma_{1}}V)$ extends  ${\sf KP}\ell^{r}$ 
with an axiom stating that
`there exists a transitive set $M$ such that $M\prec_{\Sigma_{1}}V$'.
$\Omega=\omega_{1}^{CK}$ and $\psi_{\Omega}$ is a collapsing function such that $\psi_{\Omega}(\alpha)<\Omega$.
$\mathbb{S}$ is an ordinal term denoting a stable ordinal, and 
$\Omega_{\mathbb{S}+N}$ the $N$-th admissible ordinal above $\mathbb{S}$ in the theorems.

\begin{theorem}\label{thm:2}
Suppose ${\sf KP}\ell^{r}+(M\prec_{\Sigma_{1}}V)\vdash\theta^{L_{\Omega}}$
for a $\Sigma_{1}$-sentence $\theta$.
Then we can find $n,N<\omega$ such that for $\alpha=\psi_{\Omega}(\omega_{n}(\Omega_{\mathbb{S}+N}+1))$,
$L_{\alpha}\models\theta$.
\end{theorem}

It is easy to see that a second order arithmetic 
$\Sigma^{1-}_{2}\mbox{{\rm -CA}}+\Pi^{1}_{1}\mbox{{\rm -CA}}_{0}$
 is interpreted canonically to the set theory
${\sf KP}\ell^{r}+(M\prec_{\Sigma_{1}}V)$, cf.\,subsection \ref{sec:9}.

$OT$ denotes a computable notation system of ordinals,
and $OT_{N}$ a restriction of $OT$ such that $OT=\bigcup_{0<N<\omega}OT_{N}$ and
$\psi_{\Omega}(\varepsilon_{\Omega_{\mathbb{S}+N}+1})$ denotes
the order type of $OT_{N}\cap\Omega$.
The aim of this paper is to show the following theorem, thereby the bound in Theorem \ref{thm:2}
is seen to be tight.

\begin{theorem}\label{th:wf}
For {\rm each} $N$, $\Sigma^{1-}_{2}\mbox{{\rm -CA}}+\Pi^{1}_{1}\mbox{{\rm -CA}}_{0}$
 proves that $(OT_{N},<)$ is a well ordering.
\end{theorem}

Thus the ordinal $\psi_{\Omega}(\Omega_{\mathbb{S}+\omega}):=\sup\{\psi_{\Omega}(\varepsilon_{\Omega_{\mathbb{S}+N}+1}): 0<N<\omega\}$ 
is the proof-theoretic ordinal of $\Sigma^{1-}_{2}\mbox{-CA}+\Pi^{1}_{1}\mbox{-CA}_{0}$ and of
${\sf KP}\ell^{r}+(M\prec_{\Sigma_{1}}V)$, where
$|{\sf KP}\ell^{r}+(M\prec_{\Sigma_{1}}V)|_{\Sigma_{1}^{\Omega}}$ denotes the 
$\Sigma_{1}^{\Omega}$-ordinal of ${\sf KP}\ell^{r}+(M\prec_{\Sigma_{1}}V)$, i.e., the ordinal
$\min\{\alpha\leq\omega_{1}^{CK}: \forall \theta\in\Sigma_{1}
\left(
{\sf KP}\ell^{r}+(M\prec_{\Sigma_{1}}V)\vdash\theta^{L_{\Omega}} \Rarw L_{\alpha}\models\theta
\right)\}$.

\begin{theorem}\label{th:main}
$\psi_{\Omega}(\Omega_{\mathbb{S}+\omega})  :=  \sup\{\psi_{\Omega}(\varepsilon_{\Omega_{\mathbb{S}+N}+1}): 0<N<\omega\}
=|\Sigma^{1-}_{2}\mbox{{\rm -CA}}+\Pi^{1}_{1}\mbox{{\rm -CA}}_{0}|
= |{\sf KP}\ell^{r}+(M\prec_{\Sigma_{1}}V)|_{\Sigma_{1}^{\Omega}}.
$
\end{theorem}

To prove the well-foundedness of a computable notation system,
we utilize
the distinguished class introduced by W. Buchholz\cite{Buchholz75}.

Let $OT$ be a computable notation system of ordinals with an ordinal term $\Omega_{1}$
to denote the least recursively regular ordinal $\omega_{1}^{CK}$.
Assuming that
the well-founded part $W(OT)$ of $OT$ exists as a set,
we can prove the well-foundedness of such a notation system $OT$.
When the next recursively regular ordinal $\Omega_{2}$ is in $OT$,
we further assume that 
a well-founded part $W(\mathcal{C}^{\Omega_{1}}(W_{0}))$ of a set $\mathcal{C}^{\Omega_{1}}(W_{0})$ exists,
where $W_{0}=W(OT)\cap\Omega_{1}$, and
$\alpha\in \mathcal{C}^{\Omega_{1}}(W_{0})$ iff
each component$<\!\!\Omega_{1}$ of $\alpha$ is in $W_{0}$.
Likewise when $OT$ contains $\alpha$-many terms denoting increasing sequence of
recursively regular ordinals,
we need to iterate the process of defining the well-founded parts $\alpha$-times.

Let us consider the case when there are $\alpha$-many recursively regular ordinals in $OT$
with the order type $\alpha$ of $OT$.
For example this is the case when $OT$ is a notation system for recursively inaccessible universes.
In this case the whole process should be internalized.
In other words we need a characterization of sets arising in the process.
Then \textit{distinguished sets} emerge, cf.\,(\ref{eq:distinguishedclass}) of Definition \ref{df:3wfdtg32}.
$W_{0}$ is the smallest distinguished set, and $W_{1}=W(\mathcal{C}^{\Omega_{1}}(W_{0}))\cap\Omega_{2}$
 is the next one.
Given two distinguished sets, it turns out that one is an initial segment of the other,
and the union $\mathcal{W}$ of all distinguished sets is distinguished, the \textit{maximal distinguished class}.
The maximal distinguished class $\mathcal{W}$ is $\Sigma^{1-}_{2}$-definable, and a proper class
without assuming the $\Sigma^{1}_{2}$-Comprehension Axiom.

In \cite{odMahlo, odpi3, Wienpi3d, WFnonmon2, KPpiNwfprf}, the well-foundedness of computable
notation systems is proved.
These notation systems are for recursively Mahlo universes in \cite{odMahlo},
for $\Pi_{3}$-reflecting universes in \cite{odpi3, Wienpi3d}, and
for first-order reflecting universes in \cite{WFnonmon2, KPpiNwfprf}.
In these proofs we cannot assume that the $\Sigma^{1-}_{2}$-definable class $\mathcal{W}$ exists as a set.
Instead, we introduce the maximal distinguished classes $\mathcal{W}^{P}$ inside limit universes $P$,
where $P$ is a limit of recursively regular universes, and
$\mathcal{W}^{P}$ denotes the union of all distinguished sets in $P$.
For any admissible set $Q$ with $P\in Q$, we see $\mathcal{W}^{P}\in Q$, and
rudimentary facts on distinguished sets hold on any limit universes.
Thus we can prove the well-foundedness of notation systems without assuming
the existence of the class $\mathcal{W}$ as a set.
The price to pay is an involved and intricate argument for relativized maximal classes $\mathcal{W}^{P}$.
However the maximal distinguished class $\mathcal{W}$
exists as a set in $\Sigma^{1-}_{2}\mbox{{\rm -CA}}+\Pi^{1}_{1}\mbox{{\rm -CA}}_{0}$.
For this, a crucial Lemma \ref{th:3wf16} is easier to state, and
it is not hard to see the fact $\mathbb{S}\in\mathcal{W}$ in Lemma \ref{lem:mS},
where $\mathbb{S}$ is the ordinal term denoting a stable ordinal.

\subsection{Parameter free $\Sigma^{1}_{2}${\rm -CA}}\label{sec:9}

In this subsection a second order arithmetic 
$\Sigma^{1-}_{2}\mbox{{\rm -CA}}+\Pi^{1}_{1}\mbox{{\rm -CA}}_{0}$
 is interpreted canonically to the set theory
${\sf KP}\ell^{r}+(M\prec_{\Sigma_{1}}V)$.

For subsystems of second order arithmetic, we follow largely Simpson's monograph\cite{Simpson},
and for extensions of Kripke-Platek set theory, see J\"ager's monograph\cite{J3}.
The subsystem $\Pi^{1}_{1}\mbox{{\rm -CA}}_{0}$ for $\Pi^{1}_{1}$-Comprehension axiom with restricted induction
 is the strongest one in the big five.

Let $\Sigma^{1-}_{2}\mbox{{\rm -CA}}+\Pi^{1}_{1}\mbox{{\rm -CA}}_{0}$
 denote a second order arithmetic obtained from $\Pi^{1}_{1}\mbox{{\rm -CA}}_{0}$
 by adding the axiom $\Sigma^{1-}_{2}\mbox{{\rm -CA}}$ for each parameter free\footnote{This means that no second order free variable occurs in $A$. First order parameters may occur in it.} $\Sigma^{1}_{2}$-formula $A(x)$, namely the axiom
 $\exists Y\forall x[x\in Y\lrarw A(x)]$.

The set theory ${\sf KP}\ell^{r}$ in \cite{J3} is obtained from the Kripke-Platek set theory ${\sf KP}\omega$ with
the axiom of Infinity, cf.\cite{Ba, J3}, 
by deleting $\Delta_{0}$-Collection on the universe,
restricting Foundation schema to $\Delta_{0}$-formulas, and
adding an axiom $(Lim)$, $\forall x\exists y(x\in y\land Ad(y))$, stating that
the universe is a limit of admissible sets, where $Ad$ is a unary predicate such that $Ad(y)$ is intended to designate that
`$y$ is a (transitive and) admissible set'.
Then a set theory ${\sf KP}\ell^{r}+(M\prec_{\Sigma_{1}}V)$ extends ${\sf KP}\ell^{r}$
 by adding an individual constant $M$ and
the axioms for the constant $M$:
$M$ is non-empty $M\neq\emptyset$, transitive $\forall x\in M\forall y\in x(y\in M)$, and stable $M\prec_{\Sigma_{1}}V$ for the universe $V$. $M\prec_{\Sigma_{1}}V$ means that
\begin{equation}\label{eq:stable}
 \varphi(u_{1},\ldots,u_{n}) \land \{u_{1},\ldots,u_{n}\}\subset M \to \varphi^{M}(u_{1},\ldots,u_{n})
\end{equation}
for each $\Sigma_{1}$-formula $\varphi$ in the set-theoretic language.

Note that $M$ is a model of ${\sf KP}\omega$ and the axiom $(Lim)$, i.e., a model of ${\sf KP}i$.

For a formula $A$ in the language of second order arithmetic let $A^{set}$ denote the formula obtained from $A$ by interpreting the first order variable $x$ as $x\in \omega$ and the second order variable $X$ as $X\subset\omega$.

\begin{lemma}\label{lem:9Quant}{\rm (Quantifier Theorem in p.125 of \cite{J3})}\\
For each $\Sigma^{1-}_{2}$-formula $F(n)$ without set parameters,
there exists a $\Sigma_{1}$-formula $A_{\Sigma}(n)$ in the language of set theory such that for
$F_{\Sigma}(n) :\Lrarw \exists d[Ad(d) \land A^{d}_{\Sigma}(n)]$,
\[
{\sf KP}l^{r}\vdash n\in\omega \to \{F^{set}(n)\lrarw F_{\Sigma}(n)\} 
.\]
\end{lemma}

\begin{lemma}\label{lem:9scdset}
For each sentence $A$ in the language of second order arithmetic,
\[
\Sigma^{1-}_{2}\mbox{{\rm {\rm -CA}}}+\Pi^{1}_{1}\mbox{{\rm {\rm -CA}}}_{0}\vdash A\Rarw 
{\sf KP}\ell^{r}+(M\prec_{\Sigma_{1}}V)\vdash A^{set}
.\]
\end{lemma}
{\it Proof}.\hspace{2mm}
It suffices to show that the class $\{n\in\omega: F^{set}(n)\}$ exists as a set for parameter free $\Sigma^{1-}_{2}$-formulas $F$.

By the Quantifier Theorem \ref{lem:9Quant}
pick a $\Sigma_{1}$-formula $F_{\Sigma}(n)$  such that
for $n\in\omega$,
$F^{set}(n) \lrarw F_{\Sigma}(n)$.
Since $\omega\subset M$, we obtain by the axiom (\ref{eq:stable}),
$F_{\Sigma}(n)\lrarw F_{\Sigma}^{M}(n)$ for any $n\in\omega$.
In other words $\{n\in\omega: F^{set}(n)\}=\{n\in\omega: F_{\Sigma}^{M}(n)\}$.
$\Delta_{0}$-Separation yields $\exists x(x=\{n\in\omega: F^{set}(n)\})$.
\hspace*{\fill} $\Box$
\\

Let us mention the contents of this paper.
In the next section \ref{sect:ordinalnotation} let us reproduce from \cite{singlestable}
a simultaneous definition of
iterated Skolem hulls $\mathcal{H}_{\alpha}(X)$ of sets $X$ of ordinals, ordinals 
$\psi^{f}_{\kappa}(\alpha)$ for regular cardinals 
$\kappa$, and finite functions $f$,
and Mahlo classes $Mh^{\alpha}_{k}(\xi)$.
In subsection \ref{subsec:decidable} a computable notation system $OT_{N}$ of ordinals is defined
for each positive integer $N$.
The well-foundedness of the notation system $OT_{N}$ is proved in
section \ref{sec:distinguished}.
\\

IH denotes the Induction Hypothesis, MIH the Main IH and SIH the Subsidiary IH.
Throughout this paper $N$ denotes a fixed positive integer.

\section{Ordinals for one stable ordinal}\label{sect:ordinalnotation}
For ordinals $\alpha\geq\beta$,
$\alpha-\beta$ denotes the ordinal $\gamma$ such that $\alpha=\beta+\gamma$.
Let $\alpha$ and $\beta$ be ordinals.
$\alpha\dot{+}\beta$ denotes the sum $\alpha+\beta$
when $\alpha+\beta$ equals to the commutative (natural) sum $\alpha\#\beta$, i.e., when
either $\alpha=0$ or $\alpha=\alpha_{0}+\omega^{\alpha_{1}}$ with
$\omega^{\alpha_{1}+1}>\beta$.

$\mathbb{S}$ denotes a weakly inaccessible cardinal with $\omega_{\mathbb{S}}=\mathbb{S}$, and
$\omega_{\mathbb{S}+n}$ the $n$-th regular cardinal above $\mathbb{S}$ for $0<n<\omega$.
Let $\mathbb{K}=\omega_{\mathbb{S}+N}$ for a fixed positive integer $N$.

\begin{definition}\label{df:Lam}
{\rm
For a positive integer $N$ let $\Lambda=\mathbb{K}=\omega_{\mathbb{S}+N}$.
$\varphi_{b}(\xi)$ denotes the binary Veblen function on 
$\Lambda^{+}=\mathbb{K}^{+}=\omega_{\mathbb{S}+N+1}$ with $\varphi_{0}(\xi)=\omega^{\xi}$, and
$\tilde{\varphi}_{b}(\xi):=\varphi_{b}(\Lambda\cdot \xi)$ for the epsilon number 
$\Lambda$.
Each ordinal $\varphi_{b}(\xi)$ is a fixed point of the function $\varphi_{c}$ for $c<b$
in the sense that $\varphi_{c}(\varphi_{b}(\xi))=\varphi_{b}(\xi)$.
The same holds for $\tilde{\varphi}_{b}$ and $\tilde{\varphi}_{c}$ with $c<b$.

Let $b,\xi<\Lambda^{+}$.
$\theta_{b}(\xi)$ [$\tilde{\theta}_{b}(\xi)$] denotes
a $b$-th iterate of $\varphi_{0}(\xi)=\omega^{\xi}$ [of $\tilde{\varphi}_{0}(\xi)=\Lambda^{\xi}$], resp.
Specifically ordinals
$\theta_{b}(\xi), \tilde{\theta}_{b}(\xi)<\Lambda^{+}$ are defined by recursion on $b$ as follows.
$\theta_{0}(\xi)=\tilde{\theta}_{0}(\xi)=\xi$, $\theta_{\omega^{b}}(\xi)=\varphi_{b}(\xi)$,
$\tilde{\theta}_{\omega^{b}}(\xi)=\tilde{\varphi}_{b}(\xi)$, and
$\theta_{c\dot{+}\omega^{b}}(\xi)=\theta_{c}(\theta_{\omega^{b}}(\xi))$
$\tilde{\theta}_{c\dot{+}\omega^{b}}(\xi)=\tilde{\theta}_{c}(\tilde{\theta}_{\omega^{b}}(\xi))$.

$\alpha>0$ is an additive principal number 
if $\forall \beta,\gamma<\alpha(\beta+\gamma<\alpha)$.
$\alpha>0$ is a strongly critical number if $\forall b,\xi<\alpha(\varphi_{b}(\xi)<\alpha)$.
$\Gamma_{a}$ denotes the $a$-th strongly critical number.
}
\end{definition}

It is easy to see that $\tilde{\theta}_{b+c}(\xi)=\tilde{\theta}_{b}(\tilde{\theta}_{c}(\xi))$.

Let us define a normal form of non-zero ordinals $\xi<\varphi_{\Lambda}(0)$.
Let $\xi=\Lambda^{\zeta}$.
If $\zeta<\Lambda^{\zeta}$, then $\tilde{\theta}_{1}(\zeta)$ is the normal form of $\xi$,
denoted by
$\xi=_{NF}\tilde{\theta}_{1}(\zeta)$.
Assume $\zeta=\Lambda^{\zeta}$, and let $b>0$ be the maximal ordinal such that
there exists an ordinal $\eta$ with
$\zeta=\tilde{\varphi}_{b}(\eta)>\eta$.
Then $\xi=\tilde{\varphi}_{b}(\eta)=_{NF}\tilde{\theta}_{\omega^{b}}(\eta)$.

Let $\xi=\Lambda^{\zeta_{m}}a_{m}+\cdots+\Lambda^{\zeta_{0}}a_{0}$, where
$\zeta_{m}>\cdots>\zeta_{0}$ and $0<a_{0},\ldots,a_{m}<\Lambda$.
Let $\Lambda^{\zeta_{i}}=_{NF}\tilde{\theta}_{b_{i}}(\eta_{i})$ with
$b_{i}=\omega^{c_{i}}$ for each $i$.
Then
$\xi=_{NF}\tilde{\theta}_{b_{m}}(\eta_{m})\cdot a_{m}+\cdots+\tilde{\theta}_{b_{0}}(\eta_{0})\cdot a_{0}$.

\begin{definition}\label{df:Lam2}
{\rm
Let $\xi<\varphi_{\Lambda}(0)$ be a non-zero ordinal with its normal form:
\begin{equation}\label{eq:CantornfLam}
\xi=\sum_{i\leq m}\tilde{\theta}_{b_{i}}(\xi_{i})\cdot a_{i}=_{NF}
\tilde{\theta}_{b_{m}}(\xi_{m})\cdot a_{m}+\cdots+\tilde{\theta}_{b_{0}}(\xi_{0})\cdot a_{0}
\end{equation}
where
$\tilde{\theta}_{b_{i}}(\xi_{i})>\xi_{i}$,
$\tilde{\theta}_{b_{m}}(\xi_{m})>\cdots>\tilde{\theta}_{b_{0}}(\xi_{0})$, 
$b_{i}=\omega^{c_{i}}<\Lambda$, and
$0<a_{0},\ldots,a_{m}<\Lambda$.
$\tilde{\theta}_{b_{0}}(\xi_{0})$ is said to be the \textit{tail} of $\xi$, denoted 
$\tilde{\theta}_{b_{0}}(\xi_{0})=tl(\xi)$, and
$\tilde{\theta}_{b_{m}}(\xi_{m})$ the \textit{head} of $\xi$, denoted 
$\tilde{\theta}_{b_{m}}(\xi_{m})=hd(\xi)$.

\begin{enumerate}
\item\label{df:Exp2.3}
 $\zeta$ is a \textit{part} of $\xi$
 iff there exists an $n\, (0\leq n\leq m+1)$
 such that
 $\zeta=_{NF}\sum_{i\geq n}\tilde{\theta}_{b_{i}}(\xi_{i})=
 \tilde{\theta}_{b_{m}}(\xi_{m})+\cdots+\tilde{\theta}_{b_{n}}(\xi_{n})$
 for $\xi$ in (\ref{eq:CantornfLam}).

\item\label{df:thtminus}
Let $\zeta=_{NF}\tilde{\theta}_{b}(\xi)$ with $\tilde{\theta}_{b}(\xi)>\xi$ and $b=\omega^{b_{0}}$,
and $c$ be ordinals.
An ordinal $\tilde{\theta}_{-c}(\zeta)$ is defined recursively as follows.
If $b\geq c$, then $\tilde{\theta}_{-c}(\zeta)=\tilde{\theta}_{b-c}(\xi)$.
Let $c> b$.
If $\xi>0$, then
$\tilde{\theta}_{-c}(\zeta)=\tilde{\theta}_{-(c-b)}(\tilde{\theta}_{b_{m}}(\xi_{m}))$ for the head term 
$hd(\xi)=\tilde{\theta}_{b_{m}}(\xi_{m})$ of 
$\xi$ in (\ref{eq:CantornfLam}).
If $\xi=0$, then let $\tilde{\theta}_{-c}(\zeta)=0$.

\end{enumerate}
}
\end{definition}

\begin{proposition}\label{prp:tht-}
Let $\zeta=_{NF}\tilde{\theta}_{b_{0}}(\zeta_{0})>\zeta_{0}$
and $b_{0}=\omega^{c_{0}}$.
\benu
\item\label{prp:tht41}
$\tilde{\theta}_{-c}(\zeta)
\leq\zeta$.

\item\label{prp:tht4}
Let $\xi=_{NF}\tilde{\theta}_{b_{1}}(\xi_{0})$. Then
$\zeta\leq\xi \Rarw \tilde{\theta}_{-c}(\zeta)\leq\tilde{\theta}_{-c}(\xi)$.

\item\label{prp:tht6}
$\tilde{\theta}_{c}(\tilde{\theta}_{-c}(\zeta))\leq\zeta$.
For $c\leq b_{0}$,
$\zeta=\tilde{\theta}_{c}(\tilde{\theta}_{-c}(\zeta))$ holds.

\item\label{prp:tht65}
If $c>0$, then
$\tilde{\theta}_{-(b+c)}(\zeta)=\tilde{\theta}_{-c}(hd(\tilde{\theta}_{-b}(\zeta)))$.

\item\label{prp:zigzag1}
Assume
$\tilde{\theta}_{d}(\eta)< \tilde{\theta}_{b_{0}}(\zeta_{1})$ for
$\eta>0$. Then $\tilde{\theta}_{-d}(\tilde{\theta}_{b_{0}}(\zeta_{1}))$ 
is an additive principal number such that
$\eta<\tilde{\theta}_{-d}(\tilde{\theta}_{b_{0}}(\zeta_{1}))$.
\eenu
\end{proposition}
{\it Proof}.\hspace{2mm}
\ref{prp:tht-}.\ref{prp:tht41} is seen by induction on $\zeta$.
\\

\noindent
\ref{prp:tht-}.\ref{prp:tht4} by induction on the sum $\xi+\zeta$.
Let $\xi=_{NF}\tilde{\theta}_{b_{1}}(\xi_{0})$
with $\tilde{\theta}_{b_{1}}(\xi_{0})>\xi_{0}$ and $b_{1}=\omega^{c_{1}}$, and
$\zeta<\xi$. We show $\tilde{\theta}_{-c}(\zeta)\leq\tilde{\theta}_{-c}(\xi)$.
First assume $b=b_{0}=b_{1}$. 
Then $\zeta_{0}<\xi_{0}$ and
$\tilde{\theta}_{-c}(\zeta)=\tilde{\theta}_{b-c}(\zeta_{0})\leq\tilde{\theta}_{b-c}(\xi_{0})=\tilde{\theta}_{-c}(\xi)$ when $b\geq c$.
If $b<c$, then $hd(\zeta_{0})\leq hd(\xi_{0})$ yields
$\tilde{\theta}_{-c}(\zeta)=\tilde{\theta}_{-(c-b)}(hd(\zeta_{0}))\leq\tilde{\theta}_{-(c-b)}(hd(\xi_{0}))=\tilde{\theta}_{-c}(\xi)$ by IH.

Second let $b_{0}<b_{1}$. 
We obtain $\zeta_{0}<\xi$.
If $b_{1}> c$, then
$\tilde{\theta}_{-c}(\zeta)\leq\zeta<\xi=\tilde{\theta}_{-c}(\xi)$ by Proposition \ref{prp:tht-}.\ref{prp:tht41}.
Let $b_{0}+c=c\geq b_{1}>b_{0}$.
We can assume $\zeta_{0}>0$.
We obtain $hd(\zeta_{0})\leq\zeta_{0}<\xi$, and
$\tilde{\theta}_{-c}(\zeta)=\tilde{\theta}_{-c}(hd(\zeta_{0}))\leq\tilde{\theta}_{-c}(\xi)$ by IH.

Third assume $b_{0}>b_{1}$.
We obtain $\zeta<\xi_{0}$, and
$\zeta\leq hd(\xi_{0})$.
If $c\leq b_{1}$, then 
$\tilde{\theta}_{-c}(\zeta)=\zeta<\xi_{0}\leq\tilde{\theta}_{b_{1}-c}(\xi_{0})=\tilde{\theta}_{-c}(\xi)$.
Let $b_{0}>c>b_{1}$.
Then 
$\tilde{\theta}_{-c}(\zeta)=\zeta=\tilde{\theta}_{-(c-b_{1})}(\zeta)\leq
\tilde{\theta}_{-(c-b_{1})}(hd(\xi_{0}))=
\tilde{\theta}_{-c}(\xi)$
by IH.
Finally let $c=b_{1}+c\geq b_{0}$.
Then 
$\tilde{\theta}_{-c}(\zeta)\leq \tilde{\theta}_{-c}(hd(\xi_{0}))=\tilde{\theta}_{-c}(\xi)$ by IH.
\\

\noindent
\ref{prp:tht-}.\ref{prp:tht6} by induction on $\zeta$.
If $c=b_{0}$, then
$\tilde{\theta}_{c}(\tilde{\theta}_{-c}(\zeta))=
\tilde{\theta}_{c}(\zeta_{0})=\zeta$.
If $c< b_{0}$, then 
$\tilde{\theta}_{c}(\tilde{\theta}_{-c}(\zeta))=\tilde{\theta}_{c}(\zeta)=\zeta$.
Let $c=b_{0}+e>b_{0}$. Then $\tilde{\theta}_{-c}(\zeta)=\tilde{\theta}_{-e}(hd(\zeta_{0}))$.
By IH we obtain
$\tilde{\theta}_{e}(\tilde{\theta}_{-e}(hd(\zeta_{0})))\leq hd(\zeta_{0})\leq\zeta_{0}$, and
$\tilde{\theta}_{c}(\tilde{\theta}_{-c}(\zeta))=\tilde{\theta}_{c}(\tilde{\theta}_{-e}(hd(\zeta_{0})))\leq
\tilde{\theta}_{b_{0}}(hd(\zeta_{0}))\leq\zeta$.
\\

\noindent
\ref{prp:tht-}.\ref{prp:tht65}. 
Let $c>0$.
If $b+c<b_{0}$, then
$\tilde{\theta}_{-(b+c)}(\zeta)=\zeta=\tilde{\theta}_{-c}(\tilde{\theta}_{-b}(\zeta))=
\tilde{\theta}_{-c}(hd(\tilde{\theta}_{-b}(\zeta)))$ with
$hd(\zeta)=\zeta$.
If $b_{0}=b+c=c$, then 
$\tilde{\theta}_{-(b+c)}(\zeta)=\zeta_{0}=\tilde{\theta}_{-c}(\zeta)=\tilde{\theta}_{-c}(hd(\tilde{\theta}_{-b}(\zeta)))$.
Let $b_{0}< b+c$.
If $b=b_{0}$, then
$\tilde{\theta}_{-(b+c)}(\zeta)=\tilde{\theta}_{-c}(hd(\zeta_{0}))=\tilde{\theta}_{-c}(hd(\tilde{\theta}_{-b}(\zeta)))$.
Let $b=b_{0}+e>b_{0}$.
Then
$\tilde{\theta}_{-(b+c)}(\zeta)=\tilde{\theta}_{-(e+c)}(hd(\zeta_{0}))=
\tilde{\theta}_{-c}(hd(\tilde{\theta}_{-e}(hd(\zeta_{0}))))=
\tilde{\theta}_{-c}(hd(\tilde{\theta}_{-b}(\zeta)))$ by IH.
Finally let $b<b_{0}< b+c=c$.
Then $\tilde{\theta}_{-(b+c)}(\zeta)=\tilde{\theta}_{-c}(\zeta)=
\tilde{\theta}_{-c}(hd(\tilde{\theta}_{-b}(\zeta)))$ by $\zeta=hd(\tilde{\theta}_{-b}(\zeta))$.
\\

\noindent
\ref{prp:tht-}.\ref{prp:zigzag1}.
Let us define ordinals $\xi_{i},b_{i}\,(i\leq k_{0})$ by
$\xi_{0}=\tilde{\theta}_{b_{0}}(\zeta_{1})$, and
$\xi_{i+1}=hd(\zeta_{i+1})$ for
$\xi_{i}=_{NF}\tilde{\theta}_{b_{i}}(\zeta_{i+1})$
if $\zeta_{i+1}>0$.
Otherwise $\xi_{i+1}$ is undefined, and $k_{0}=i$.
We claim that $b_{0}+\cdots+b_{k_{0}}>d$.
Suppose $b=b_{0}+\cdots+b_{k_{0}}\leq d$, and let $b_{k_{0}+1}=d-b$ and
$\eta_{i}=\tilde{\theta}_{b_{i}+\cdots+b_{k_{0}}+b_{k_{0}+1}}(\eta)$.
Then the assumption yields
$\eta_{0}=\tilde{\theta}_{b_{0}}(\eta_{1})=
\tilde{\theta}_{d}(\eta)\leq\xi_{0}=\tilde{\theta}_{b_{0}}(\zeta_{1})$.
We obtain 
$\eta_{1}\leq\zeta_{1}$, and
$\eta_{1}=\tilde{\theta}_{b_{1}}(\eta_{2})
\leq hd(\zeta_{1})=\xi_{1}=\tilde{\theta}_{b_{1}}(\zeta_{2})$
for the additive principal number $\eta_{1}$.
We see inductively that 
$\tilde{\theta}_{b_{k_{0}}}(\tilde{\theta}_{b_{k_{0}+1}}(\eta))=\eta_{k_{0}}\leq
\xi_{k_{0}}=\tilde{\theta}_{b_{k_{0}}}(\zeta_{k_{0}+1})$, and
$0<\tilde{\theta}_{b_{k_{0}+1}}(\eta)\leq\zeta_{k_{0}+1}=0$. This is a contradiction, and $b_{0}+\cdots+b_{k_{0}}>d$.
Therefore $\tilde{\theta}_{-d}(\tilde{\theta}_{b_{0}}(\zeta_{1}))$ 
is an additive principal number.

For the least number $k=\min\{k\leq k_{0}: b_{0}+\cdots+b_{k}>d\}$,
let
$\{e_{1}<\cdots <e_{n}\}=\{b_{0}+\cdots+b_{i}:i<k\}\cup\{d\}\, (k\leq n\leq k+1)$, 
$e_{0}=0$
and $d_{i}=e_{i}-e_{i-1}$.
Then $0<d_{i}\leq b_{i-1}$.

Now suppose $\eta\geq\tilde{\theta}_{-d}(\tilde{\theta}_{b_{0}}(\zeta_{1}))$.
By Propositions \ref{prp:tht-}.\ref{prp:tht6} and \ref{prp:tht-}.\ref{prp:tht65} 
we obtain
$\tilde{\theta}_{d_{n}}(\eta)\geq\tilde{\theta}_{d_{n}}(\tilde{\theta}_{-d}(\tilde{\theta}_{b_{0}}(\zeta_{1})))
=\tilde{\theta}_{d_{n}}(\tilde{\theta}_{-d_{n}}(hd(\tilde{\theta}_{-e_{n-1}}(\tilde{\theta}_{b_{0}}(\zeta_{1})))))
\geq\tilde{\theta}_{-e_{n-1}}(\tilde{\theta}_{b_{0}}(\zeta_{1}))$, where
$d=e_{n-1}+d_{n}$, $0<d_{n}\leq b_{n-1}$ and
$\tilde{\theta}_{-e_{n-1}}(\tilde{\theta}_{b_{0}}(\zeta_{1}))=\tilde{\theta}_{b_{n-1}}(\zeta_{n})$.
We see inductively that
$\tilde{\theta}_{d-e_{i-1}}(\eta)\geq\tilde{\theta}_{-e_{i-1}}(\tilde{\theta}_{b_{0}}(\zeta_{1}))$
for $d-e_{i-1}=d_{i}+\cdots+d_{n}$.
We arrive at a contradiction
$\tilde{\theta}_{d}(\eta)\geq\tilde{\theta}_{-e_{0}}(\tilde{\theta}_{b_{0}}(\zeta_{1}))=\tilde{\theta}_{b_{0}}(\zeta_{1})$ for $e_{0}=0$.
\hspace*{\fill} $\Box$
\\

A `Mahlo degree' $m(\pi)$ of ordinals $\pi$ is defined to be a finite function 
$f:\mathbb{K} \to\varphi_{\mathbb{K}}(0)$.
$a,b,c,\alpha,\beta,\gamma,\delta,\ldots$ range over ordinals$<\!\varepsilon_{\mathbb{K}+1}$,
$\xi,\zeta,\eta,\ldots$ range over ordinals$<\!\! \varphi_{\mathbb{K}}(0)$.

\begin{definition}\label{df:Lam3}
{\rm

  \benu
  \item
A function $f:\mathbb{K}\to \varphi_{\mathbb{K}}(0)$ with a \textit{finite} support
${\rm supp}(f)=\{c<\mathbb{K}: f(c)\neq 0\}\subset \mathbb{K}$ is said to be a \textit{finite function}
if
$\forall i>0(a_{i}=1)$ and $a_{0}=1$ when $b_{0}>1$
in
$f(c)=_{NF}\tilde{\theta}_{b_{m}}(\xi_{m})\cdot a_{m}+\cdots+\tilde{\theta}_{b_{0}}(\xi_{0})\cdot a_{0}$
for any $c\in{\rm supp}(f)$.
A function $f$ with a finite support is identified with the finite function $f\restrict \supp(f)$.
When $c\not\in {\rm supp}(f)$, let $f(c):=0$.
$K(f):=\bigcup\{\{c,f(c)\}: c\in \supp(f)\}$.
$f,g,h,\ldots$ range over finite functions.

For an ordinal $c$, $f_{c}$ and $f^{c}$ are restrictions of $f$ to the domains
${\rm supp}(f_{c})=\{d\in\supp(f): d< c\}$ and ${\rm supp}(f^{c})=\{d\in\supp(f): d\geq c\}$.
$g_{c}*f^{c}$ denotes the concatenated function such that
${\rm supp}(g_{c}*f^{c})={\rm supp}(g_{c})\cup {\rm supp}(f^{c})$, $(g_{c}*f^{c})(a)=g(a)$ for $a<c$, and
$(g_{c}*f^{c})(a)=f(a)$ for $a\geq c$.

\item\label{df:Exp2.5}
Let $f$ be a finite function and $c,\xi$ ordinals.
A relation $f<^{c}\xi$ is defined by induction on the
cardinality of the finite set $\{d\in {\rm supp}(f): d>c\}$ as follows.
If $f^{c}=\emptyset$, then $f<^{c}\xi$ holds.
For $f^{c}\neq\emptyset$,
$f<^{c}\xi$ iff
there exists a part $\mu$ of $\xi$ such that
$f(c)< \mu$
and 
$f<^{c+d}\tilde{\theta}_{-d}(tl(\mu))$ 
for $d=\min\{c+d\in \supp(f): d>0\}$.

\eenu

}
\end{definition}
For example, let $f$ be a finite function with
${\rm supp}(f)=\{0,1\}$, and
$\xi=_{NF}\tilde{\theta}_{1}(\xi_{1})+\tilde{\theta}_{1}(\xi_{0})$.
Then $f<^{0}\xi$ iff either 
$f(0)<\xi \spand f(1)<\tilde{\theta}_{-1}(tl(\xi))=\tilde{\theta}_{-1}(\tilde{\theta}_{1}(\xi_{0}))=\xi_{0}$, or
$f(0)<\tilde{\theta}_{1}(\xi_{1}) \spand f(1)<\tilde{\theta}_{-1}(\tilde{\theta}_{1}(\xi_{1}))=\xi_{1}$.
If $f(0)=_{NF}\tilde{\theta}_{1}(\xi_{1})+\tilde{\theta}_{1}(\zeta)$ with $\zeta<\xi_{0}\leq f(1)$, then
$f\not<^{0}\xi$.

From Proposition \ref{prp:tht-}.\ref{prp:tht4} we see that, cf.\cite{singlestable}
\begin{equation}\label{prp:idless}
f<^{c}\xi\leq\zeta \Rightarrow f<^{c}\zeta
\end{equation}

\subsection{Skolem hulls and collapsing functions}
In this subsection let us reproduce from \cite{singlestable}
a simultaneous definition of
Skolem hulls $\mathcal{H}_{a}(\alpha)$, collapsing functions $\psi_{\pi}^{f}(\alpha)$ and
Mahlo classes $Mh^{a}_{c}(\xi)$.

\begin{definition}
{\rm
\benu

\item\label{df:Lam3}
Let $A\subset\mathbb{S}$ be a set, and $\alpha\leq\mathbb{S}$ a limit ordinal.
\[
\alpha\in M(A) :\Lrarw A\cap\alpha \mbox{ is stationary in }
\alpha
\Lrarw \mbox{ every club subset of } \alpha \mbox{ meets } A.
\]

 \item\label{df:Lam4}
 $\kappa^{+}$ denotes the next regular ordinal above $\kappa$.
 For $n<\omega$, $\kappa^{+n}$ is defined recursively by
 $\kappa^{+0}=\kappa$ and $\kappa^{+(n+1)}=\left(\kappa^{+n}\right)^{+}$.

\item
$\Omega_{\alpha}:=\omega_{\alpha}$ for $\alpha>0$, $\Omega_{0}:=0$, and
$\Omega=\Omega_{1}$.
$\mathbb{S}$ is a weakly inaccessible cardinal, and
$\Omega_{\mathbb{S}}=\mathbb{S}$.
$\Omega_{\mathbb{S}+n}=\mathbb{S}^{+n}$ is the $n$-th cardinal above $\mathbb{S}$.
\eenu
}
\end{definition}

In the following Definition \ref{df:Cpsiregularsm}, 
$\varphi\alpha\beta=\varphi_{\alpha}(\beta)$ denotes the binary Veblen function on $\mathbb{K}^{+}=\omega_{\mathbb{S}+N+1}$,
$\tilde{\theta}_{b}(\xi)$ the function defined in Definition \ref{df:Lam}
for $\Lambda=\mathbb{K}=\omega_{\mathbb{S}+N}$ with the positive integer $N$.

For 
$a<\varepsilon_{\mathbb{K}+1}$,
$c<\mathbb{K}$, and
$\xi<\Gamma_{\mathbb{K}+1}$, 
define simultaneously 
classes $\mathcal{H}_{a}(X)\subset\Gamma_{\mathbb{K}+1}$,
$Mh^{a}_{c}(\xi)\subset(\mathbb{S}+1)$, and 
ordinals $\psi_{\kappa}^{f}(a)\leq\kappa$ by recursion on ordinals $a$ as follows.

\begin{definition}\label{df:Cpsiregularsm}
{\rm
Let
$\mathbb{K}=\Omega_{\mathbb{S}+N}$, 
$\mathcal{H}_{a}[Y](X):=\mathcal{H}_{a}(Y\cup X)
$ for sets $Y\subset \Gamma_{\mathbb{K}+1}$.
Let $a<\varepsilon_{\mathbb{K}+1}$ and $X\subset\Gamma_{\mathbb{K}+1}$.

\begin{enumerate}
\item\label{df:Cpsiregularsm.1}
(Inductive definition of $\mathcal{H}_{a}(X)$.)

\begin{enumerate}
\item\label{df:Cpsiregularsm.10}
$\{0,\Omega_{1},\mathbb{S}\}\cup\{\Omega_{\mathbb{S}+n}: 0<n\leq N\}\cup X\subset\mathcal{H}_{a}(X)$.

\item\label{df:Cpsiregularsm.11}
If $x, y \in \mathcal{H}_{a}(X)$,
then $x+y\in \mathcal{H}_{a}(X)$,
and 
$\varphi xy\in \mathcal{H}_{a}(X)$.

\item\label{df:Cpsiregularsm.12}
Let $\alpha\in\mathcal{H}_{a}(X)\cap\mathbb{S}$. Then for each $0<k\leq N$, 
$\Omega_{\alpha+k}\in\mathcal{H}_{a}(X)$.

\item\label{df:Cpsiregularsm.1345}
Let $\alpha=\psi_{\pi}^{f}(b)$ with $\{\pi,b\}\subset\mathcal{H}_{a}(X)$, 
$\pi\in\{\mathbb{S}\}\cup\Psi_{N}$,
$b<a$, and a finite function $f$ such that
$K(f)\subset\mathcal{H}_{a}(X)\cap\mathcal{H}_{b}(\alpha)$.

Then $\alpha\in\mathcal{H}_{a}(X)\cap\Psi_{N}$.

\end{enumerate}

\item\label{df:Cpsiregularsm.2}
 (Definitions of $Mh^{a}_{c}(\xi)$ and $Mh^{a}_{c}(f)$)
\\
The classes $Mh^{a}_{c}(\xi)$ are defined for $c< \mathbb{K}$,
and ordinals $a<\varepsilon_{\mathbb{K}+1}$, $\xi<\Gamma_{\mathbb{K}+1}$.
Let $\pi$ be a regular ordinal$\leq \mathbb{S}$. Then 
by main induction on ordinals $\pi\leq\mathbb{S}$
with subsidiary induction on $c<\mathbb{K}$ 
we define $\pi\in Mh^{a}_{c}(\xi)$ iff 
$\{a,c\}\cup K(\xi)\subset\mathcal{H}_{a}(\pi)$ and
\begin{equation}\label{eq:dfMhkh}
 \forall f<^{c}\xi 
 \forall g \left(
 K(f)\cup K(g) \subset\mathcal{H}_{a}(\pi) 
 \,\&\, 
\pi\in Mh^{a}_{0}(g_{c})
 \Rightarrow \pi\in M(Mh^{a}_{0}(g_{c}*f^{c}))
 \right)
\end{equation}
where $f, g$ vary through 
finite
 functions,
and 
\begin{eqnarray*}
Mh^{a}_{c}(f)  & := & \bigcap\{Mh^{a}_{d}(f(d)): d\in {\rm supp}(f^{c})\}
\\
& = &
\bigcap\{Mh^{a}_{d}(f(d)): c\leq d\in {\rm supp}(f)\}.
\end{eqnarray*}
In particular
$Mh^{a}_{0}(g_{c})=\bigcap\{Mh^{a}_{d}(g(d)): d\in {\rm supp}(g_{c})\}
=\bigcap\{Mh^{a}_{d}(g(d)): c> d\in {\rm supp}(g)\}$.
When $f=\emptyset$ or $f^{c}=\emptyset$, let $Mh^{a}_{c}(\emptyset):=\mathbb{K}$.

\item\label{df:Cpsiregularsm.3}
 (Definition of $\psi_{\pi}^{f}(a)$)
\\
 Let $a<\varepsilon_{\mathbb{K}+1}$ 
 be an ordinal, $\pi$ a regular ordinal and
 $f$ a finite function.
Then let
{\small
\begin{equation}\label{eq:Psivec}
\psi_{\pi}^{f}(a)
 :=  \min(\{\pi\}\cup\{\kappa\in Mh^{a}_{0}(f)\cap\pi:   \mathcal{H}_{a}(\kappa)\cap\pi\subset\kappa ,
   \{\pi,a\}\cup K(f)\subset\mathcal{H}_{a}(\kappa)
\})
\end{equation}
}
For the empty function $\emptyset$,
$\psi_{\pi}(a):=\psi_{\pi}^{\emptyset}(a)$.

\item
For classes $A\subset(\mathbb{S}+1)$, let
$\alpha\in M^{a}_{c}(A)$ iff $\alpha\in A$ and
\begin{equation}\label{eq:Mca}
\forall g
[
\alpha\in Mh_{0}^{a}(g_{c}) \,\&\, K(g_{c})\subset\mathcal{H}_{a}(\alpha) \Rightarrow
\alpha\in M\left( Mh_{0}^{a}(g_{c}) \cap A \right)
]
\end{equation}
\end{enumerate}
}

\end{definition}

Assuming an existence of a 
shrewd cardinal introduced by M. Rathjen\cite{RathjenAFML2},
we show in \cite{singlestable} that
$\psi_{\mathbb{S}}^{f}(a)<\mathbb{S}$
if $\{a,c,\xi\}\subset\mathcal{H}_{a}(\mathbb{S})$ with 
$c<\mathbb{K}$, $a,\xi<\varepsilon_{\mathbb{K}+1}$,
and ${\rm supp}(f)=\{c\}$ and $f(c)=\xi$.
Moreover $\psi_{\pi}^{g}(b)<\pi$ 
provided that
$\pi\in Mh^{b}_{0}(f)$, $K(g)\cup\{\pi,b\}\subset\mathcal{H}_{b}(\pi)$, and
$g$ is a finite function defined from a finite function $f$ and ordinals $d,c$ as follows.
$d<c\in \supp(f)$ 
with
$(d,c)\cap \supp(f)=(d,c)\cap \supp(g)=\emptyset$, 
$g_{d}=f_{d}$,
$g(d)<f(d)+\tilde{\theta}_{c-d}(f(c))\cdot\omega$, and
$g<^{c}f(c)$.
Also the following Lemma \ref{lem:stepdown} is shown in \cite{singlestable}.

\begin{lemma}\label{lem:stepdown}
Assume $\mathbb{S}\geq\pi\in Mh^{a}_{d}(\xi)\cap Mh^{a}_{c}(\xi_{0})$, $\xi_{0}\neq 0$,
$d<c$, and
$\{a,c,d\}\subset\mathcal{H}_{a}(\pi)$. 
Moreover let $\tilde{\theta}_{c-d}(\xi_{0})\geq\xi_{1}\in\mathcal{H}_{a}(\pi)$
and $tl(\xi)>\xi_{1}$ when $\xi\neq 0$.
Then
$\pi\in Mh^{a}_{d}(\xi+\xi_{1})\cap M^{a}_{d}(Mh^{a}_{d}(\xi+\xi_{1}))$.
\end{lemma}

\subsection{Normal forms in ordinal notations}
In this subsection we introduce an \textit{irreducibility} of finite functions,
which is needed to define a normal form in ordinal notations.

The following Proposition \ref{prp:nfform}
is shown from Lemma \ref{lem:stepdown} in \cite{singlestable}.

\begin{proposition}\label{prp:nfform}
Let $f$ be a finite function such that
$\{a\}\cup K(f)\subset\mathcal{H}_{a}(\pi)$.
Assume $tl(f(c))\leq\tilde{\theta}_{d}(f(c+d))$ 
holds for some $c, c+d\in {\rm supp}(f)$ with 
$d>0$.
Then 
$\pi\in Mh^{a}_{0}(f) \Leftrightarrow \pi\in Mh^{a}_{0}(g)$ holds,
where
$g$ is a finite function such that
$g(c)=f(c)-tl(f(c))$ and $g(e)=f(e)$ for every $e\neq c$.
\end{proposition}

\begin{definition}\label{df:irreducible}
{\rm

An \textit{irreducibility} of finite functions $f$ is defined by induction on the cardinality
$n$ of the finite set $\supp(f)$.
If $n\leq 1$, $f$ is defined to be irreducible.
Let $n\geq 2$ and $c<c+d$ be the largest two elements in $\supp(f)$, and let $g$ be 
a finite function
such that $\supp(g)=\supp(f_{c})\cup\{c\}$, $g_{c}=f_{c}$ and
$g(c)=f(c)+\tilde{\theta}_{d}(f(c+d))$.
Then $f$ is irreducible iff 
$tl(f(c))> \tilde{\theta}_{d}(f(c+d))$ and
$g$ is irreducible.

}
\end{definition}

\begin{proposition}\label{prp:zigzag}
Let $f$ be a finite function.
Assume $f(c)<\xi$, and
$\tilde{\theta}_{d}(f(c+d))<tl(\xi)$ for every $c<c+d\in{\rm supp}(f)$.
Then
$f<^{c}\xi$ holds.

\end{proposition}
{\it Proof}.\hspace{2mm} 
Let $c<c+d\in \supp(f)$.
For $d_{0}>0$, we obtain
$\tilde{\theta}_{-(d+d_{0})}(tl(\xi))= \tilde{\theta}_{-d_{0}}(hd(\tilde{\theta}_{-d}(tl(\xi)))$
by Proposition \ref{prp:tht-}.\ref{prp:tht65}.
Hence the proposition follows from 
Proposition \ref{prp:tht-}.\ref{prp:zigzag1}.
\hspace*{\fill} $\Box$

\begin{definition}\label{df:lx}
 {\rm 
 Let  $f,g$ be irreducible functions, and $b,a$ ordinals.
 \begin{enumerate}
 \item\label{df:lxx}
Let us define a relation $f<^{b}_{lx}g$
by induction on the cardinality of the finite set
$\{e\in{\rm supp}(f)\cup{\rm supp}(g): e\geq b\}$ as follows.
$f<^{b}_{lx}g$ holds iff $f^{b}\neq g^{b}$ and
for the ordinal $c=\min\{c\geq b : f(c)\neq g(c)\}$,
one of the following conditions is met:

\begin{enumerate}

\item\label{df:lx.23}
$f(c)<g(c)$ and let $\mu$ be the shortest part of $g(c)$ such that $f(c)<\mu$.
Then for any $c<c+d\in{\rm supp}(f)$,  
if $tl(\mu)\leq\tilde{\theta}_{d}(f(c+d))$, then 
$f<_{lx}^{c+d}g$ holds.

\item\label{df:lx.24}
$f(c)>g(c)$ and let $\nu$ be the shortest part of $f(c)$ such that $\nu>g(c)$.
Then there exist a $c<c+d\in {\rm supp}(g)$ such that
$f<_{lx}^{c+d}g$ and
$tl(\nu)\leq \tilde{\theta}_{d}(g(c+d))$.

\end{enumerate}

\item
$Mh^{a}_{b}(f)\prec Mh^{a}_{b}(g)$ holds iff
\[
\forall\pi\in Mh^{a}_{b}(g)\forall b_{0}\leq b
\left(
K(f)\subset\mathcal{H}_{a}(\pi) \,\&\, \pi\in Mh^{a}_{b_{0}}(f_{b})
 \Rightarrow \pi\in M(Mh^{a}_{b_{0}}(f))
\right)
.
\]
\end{enumerate}
}
\end{definition}

\begin{proposition}\label{prp:psinucomparison}
Let $f,g$ be irreducible finite functions.
\begin{enumerate}
\item\label{prp:psinucomparison.1}
Let $f(c)>g(c)$ and $\nu$ be the shortest part of $f(c)$ such that $\nu>g(c)$.
Assume that there is an ordinal $d>0$ such that $c+d\in {\rm supp}(g)$ and
$tl(\nu)\leq \tilde{\theta}_{d}(g(c+d))$.
Then
$\forall e\leq d[c<c+e\in \supp(f) \Rightarrow f(c+e)<\tilde{\theta}_{d-e}(g(c+d))]$.

If $d$ is the least such one, then
$d=\min\{d>0: c+d\in{\rm supp}(g)\}$.

\item\label{prp:psinucomparison.2}
Let $f^{b}\neq g^{b}$. Then
either $f<_{lx}^{b}g$ or $g<_{lx}^{b}f$ holds.

\end{enumerate}
\end{proposition}
{\it Proof}.\hspace{2mm}
\ref{prp:psinucomparison}.\ref{prp:psinucomparison.1}.
Let $c<c+e\in \supp(f)$ with $e\leq d$.
For the irreducible function $f$ with $f(c)>0$, we see that 
$tl(f(c))>\tilde{\theta}_{e}(f(c+e))$.
Suppose $f(c+e)\geq \tilde{\theta}_{d-e}(g(c+d))$.
Then we would have
$tl(\nu)\geq tl(f(c))>\tilde{\theta}_{e}(f(c+e))\geq\tilde{\theta}_{e}(\tilde{\theta}_{d-e}(g(c+d)))=\tilde{\theta}_{d}(g(c+d))$.
This is a contradiction, and the first claim follows.

Suppose $c+e\in{\rm supp}(g)\cap(c,c+d)$.
Then by the minimality of $d$ we obtain
$\tilde{\theta}_{e}(g(c+e))<tl(\nu)\leq\tilde{\theta}_{d}(g(c+d))$, and hence
$g(c+e)<\tilde{\theta}_{d-e}(g(c+d))$.
On the other hand we have
$\tilde{\theta}_{d-e}(g(c+d))<tl((g(c+e))\leq g(c+e)$ by the irreducibility of $g$.
This is a contradiction, and the second claim follows.
\\

\noindent
\ref{prp:psinucomparison}.\ref{prp:psinucomparison.2} 
 by induction on the cardinality of the finite set
$\{e\in{\rm supp}(f)\cup{\rm supp}(g): e\geq b\}$.
Let $f(c)<g(c)$ with $c\geq b$,
 and let $\mu$ be the shortest part of $g(c)$ such that $f(c)<\mu$.
Assume $tl(\mu)\leq\tilde{\theta}_{d}(f(c+d))$ for a $c<c+d\in{\rm supp}(f)$.
Then we obtain $g(c+d)<f(c+d)$ 
by Proposition \ref{prp:psinucomparison}.\ref{prp:psinucomparison.1}, and hence
$f^{c+d}\neq g^{c+d}$.
By IH either $f<_{lx}^{c+d}g$ or $g<_{lx}^{c+d}f$ holds.
Therefore either $f<_{lx}^{b}g$ or $g<_{lx}^{b}f$ holds.
\hspace*{\fill} $\Box$

\begin{lemma}\label{lem:psinucomparison}
Let $f,g$ be irreducible finite functions, and $b$ an ordinal such that $f^{b}\neq g^{b}$.
If $f<^{b}_{lx}g$, then
$Mh^{a}_{b}(f)\prec Mh^{a}_{b}(g)$ holds for every ordinal $a$.
\end{lemma}
{\it Proof}.\hspace{2mm}
We show the lemma
 by induction on the cardinality of the finite set
$\{e\in{\rm supp}(f)\cup{\rm supp}(g): e\geq b\}$.
Let $f<^{b}_{lx}g$ and $c=\min\{c\geq b : f(c)\neq g(c)\}$.
For finite functions $h$ and ordinals $b\leq c$, let $h^{b}_{c}:=(h^{b})_{c}$ denote
the restriction of $h$ to the interval $[b,c)$.
Let $\pi\in Mh^{a}_{b}(g)$ with $K(f)\subset\mathcal{H}_{a}(\pi)$.
Also assume $\pi\in Mh^{a}_{b_{0}}(f_{b})=Mh^{a}_{0}(f^{b_{0}}_{b})$ with $b_{0}\leq b$.
\\
{\bf Case 1}. \ref{df:lx}(\ref{df:lx.23}) holds:
If $\tilde{\theta}_{d}(f(c+d))<tl(\mu)$ for any $c<c+d\in{\rm supp}(f)$,  then
$f<^{c}\mu\leq g(c)$ follows from Proposition \ref{prp:zigzag}, and $f<^{c}g(c)$ by (\ref{prp:idless}).
We have $\pi\in Mh^{a}_{c}(g(c))\cap Mh^{a}_{0}(f^{b_{0}}_{c})$ with $f^{b}_{c}=g^{b}_{c}$.
Hence we obtain $\pi\in M(Mh_{b_{0}}^{a}(f))$ for $f^{b_{0}}=f^{b_{0}}_{c}*f^{c}$ 
with $Mh_{b_{0}}^{a}(f)=Mh_{0}^{a}(f^{b_{0}})$ 
by $f<^{c}g(c)$ and the definition (\ref{eq:dfMhkh}).

Assume that there exists a $c<c+d\in{\rm supp}(f)$ such that $tl(\mu)\leq\tilde{\theta}_{d}(f(c+d))$.
Let $d$ be the least such one.
Then $d=\min\{d>0: c+d\in{\rm supp}(f)\}$ by Proposition \ref{prp:psinucomparison}.\ref{prp:psinucomparison.1}.
We obtain $f<_{lx}^{c+d}g$, and $Mh^{a}_{c+d}(f)\prec Mh^{a}_{c+d}(g)$ by IH.
We have $\pi\in Mh^{a}_{c}(g(c))\cap Mh^{a}_{0}(f^{b_{0}}_{c})$ with 
$f(c)<g(c)$ and $f^{b}_{c}=g^{b}_{c}$.
Moreover $f^{c+1}_{c+d}=\emptyset$.
Hence we obtain $\pi\in Mh^{a}_{0}(f^{b_{0}}_{c+d})$ by (\ref{prp:idless}).
$\pi\in Mh^{a}_{c+d}(g)$ yields $\pi\in M(Mh_{b_{0}}^{a}(f))$.
\\
{\bf Case 2}. \ref{df:lx}(\ref{df:lx.24}) holds: Let $f<_{lx}^{c+d}g$ 
and $tl(\nu)\leq \tilde{\theta}_{d}(g(c+d))$ for a $c<c+d\in \supp(g)$.
We have $\nu=g(c)+tl(\nu)$.

We have $\pi\in Mh^{a}_{c}(g(c))\cap Mh^{a}_{c+d}(g)\cap Mh^{a}_{0}(f^{b}_{c})$ 
with $f^{b}_{c}=g^{b}_{c}$.
Let $f(c)=\xi\dot{+}tl(\nu)\dot{+}\eta$ with
$\xi\leq g(c)$ and $\eta\leq tl(\nu)\cdot n$ for an $n<\omega$.
By (\ref{prp:idless}) we obtain $\pi\in Mh^{a}_{c}(\xi)$, and
$\pi\in Mh^{a}_{c}(f(c))$ by Lemma \ref{lem:stepdown}.
Hence by Proposition \ref{prp:psinucomparison}.\ref{prp:psinucomparison.1} and Lemma \ref{lem:stepdown} we obtain
$\pi\in Mh^{a}_{0}(f^{b_{0}}_{c+d})$.
On the other hand we have $Mh^{a}_{c+d}(f)\prec Mh^{a}_{c+d}(g)$ by IH with $f<_{lx}^{c+d}g$.
Therefore we obtain $\pi\in M(Mh_{b_{0}}^{a}(f))$.
\hspace*{\fill} $\Box$

\begin{proposition}\label{prp:psicomparison}
Let $f,g$ be 
irreducible finite functions, and assume that
$\psi_{\pi}^{f}(b)<\pi$ and $\psi_{\kappa}^{g}(a)<\kappa$.

Then $\psi_{\pi}^{f}(b)<\psi_{\kappa}^{g}(a)$ iff one of the following cases holds:
\benu
\item\label{prp:psicomparison.0}
$\pi\leq \psi_{\kappa}^{g}(a)$.

\item\label{prp:psicomparison.1}
$b<a$, $\psi_{\pi}^{f}(b)<\kappa$ and 
$K(f)\cup\{\pi,b\}\subset\mathcal{H}_{a}(\psi_{\kappa}^{g}(a))$.

\item\label{prp:psicomparison.2}
$b>a$ and $K(g)\cup\{\kappa,a\}\not\subset\mathcal{H}_{b}(\psi_{\pi}^{f}(b))$.

\item\label{prp:psicomparison.25}
$b=a$, $\kappa<\pi$ and $\kappa\not\in\mathcal{H}_{b}(\psi_{\pi}^{f}(b))$.

\item\label{prp:psicomparison.3}
$b=a$, $\pi=\kappa$, $K(f)\subset\mathcal{H}_{a}(\psi_{\kappa}^{g}(a))$, and
$f<^{0}_{lx}g$.

\item\label{prp:psicomparison.4}

$b=a$, $\pi=\kappa$, 
$K(g)\not\subset\mathcal{H}_{b}(\psi_{\pi}^{f}(b))$.

\eenu

\end{proposition}
{\it Proof}.\hspace{2mm}
This is seen as in Proposition 2.19 of \cite{KPPiN} using Lemma \ref{lem:psinucomparison}.
\hspace*{\fill} $\Box$
\\

We need the following Lemma \ref{lem:of} to show Proposition \ref{prp:psimS},
and Lemma \ref{lem:oflx} for Lemma \ref{lem:psiw}
in section \ref{sec:distinguished}.

\begin{definition}\label{df:Lam.of}
{\rm
\benu
\item
$a(\xi)$ denotes an ordinal defined recursively by
$a(0)=0$, and
$a(\xi)=\sum_{i\leq m}\tilde{\theta}_{b_{i}}(\omega\cdot a(\xi_{i}))$
when $\xi=_{NF}\sum_{i\leq m}\tilde{\theta}_{b_{i}}(\xi_{i})$ in (\ref{eq:CantornfLam}).

\item
For irreducible functions $f$ let us associate ordinals $o(f)<\Gamma_{\mathbb{K}+1}$ as follows.
$o(\emptyset)=0$ for the empty function $f=\emptyset$.
Let $\{0\}\cup \supp(f)=\{0=c_{0}<c_{1}<\cdots<c_{n}\}$,
$f(c_{i})=\xi_{i}<\Gamma_{\mathbb{K}+1}$ for $i>0$, and $\xi_{0}=0$.
Define ordinals $\zeta_{i}=o(f;c_{i})$ by
$\zeta_{n}=\omega\cdot a(\xi_{n})$, and $\zeta_{i}=\omega\cdot a(\xi_{i})+\tilde{\theta}_{c_{i+1}-c_{i}}(\zeta_{i+1}+1)$.
Finally let $o(f)=\zeta_{0}=o(f;c_{0})$.

\item
For $d\not\in \{0\}\cup \supp(f)$, let
$o(f;d)=0$ if $f^{d}=\emptyset$.
Otherwise $o(f;d)=\tilde{\theta}_{c-d}(o(f;c)+1)$ for
$c=\min( {\rm supp}(f^{d}))$.
\eenu
}
\end{definition}
The following Proposition \ref{prp:ofb} is easily seen.

\begin{proposition}\label{prp:ofb}
\benu
\item\label{prp:ofb1}
$\xi<\zeta \Rarw a(\xi)<a(\zeta)$.

\item\label{prp:ofb2}
$tl(a(\xi))=a(tl(\xi))$.

\item\label{prp:ofb3}
$a(\tilde{\theta}_{-c}(\xi))=\tilde{\theta}_{-c}(a(\xi))$.
\eenu
\end{proposition}

\begin{proposition}\label{prp:oflx}
Let $f$ and $g$ be irreducible finite functions, and $c,d$ ordinals.

\begin{enumerate}
\item\label{prp:oflx.05}
Let $\{d<e\}\subset{\rm supp}(g)$ with $(d,e)\cap{\rm supp}(g)=\emptyset$.
Then
$\omega\cdot tl(a(g(d)))>\tilde{\theta}_{e-d}(o(g;e)+\omega)$.

\item\label{prp:oflx.0}
Let ${\rm supp}(g)\ni d<c$ with $g^{c}\neq\emptyset$.
Then $\tilde{\theta}_{c-d}(o(g;c)+1)<o(g;d)$.

\item\label{prp:oflx.1}
If $f<^{d}g(d)$, then $o(f;d)+\omega\leq o(g;d)$.

\item\label{prp:oflx.2}
Let $f(c)<g(c)+\tilde{\theta}_{d-c}(g(d))\cdot \omega$, 
$o(f;d)+\omega\leq o(g;d)$ and
$(c,d)\cap{\rm supp}(g)=\emptyset$ for ${\rm supp}(f)\ni c<d\in{\rm supp}(g)$.
Moreover if $(c,d)\cap{\rm supp}(f)\neq\emptyset$, then
$f(c)>g(c)$ and $tl(\nu)\leq\tilde{\theta}_{d-c}(g(d))$ for the shortest part $\nu$
of $f(c)$ such that $\nu>g(c)$.
Then $o(f;c)+\omega\leq o(g;c)$.

\item\label{prp:oflx.3}
Let $f(d)<g(d)$, $g^{c}\neq\emptyset$, 
$o(f;c)+\omega\leq o(g;c)$ and
$(d,c)\cap{\rm supp}(f)=\emptyset$ for ${\rm supp}(g)\ni d<c\in{\rm supp}(f)$.
Then $o(f;d)+\omega\leq o(g;d)$.

\item\label{prp:oflx.4}
Assume $o(f;c)+\omega\leq o(g;c)$ and 
$f_{c}=g_{c}$ for $c\in{\rm supp}(f)\cup{\rm supp}(g)$.
Then $o(f)<o(g)$.

\end{enumerate}
\end{proposition}
{\it Proof}.\hspace{2mm}
\ref{prp:oflx}.\ref{prp:oflx.05}.
Let $\{d=e_{0}<e_{1}=e<\cdots<e_{n}\}={\rm supp}(g^{d})$.
By the irreducibility of $g$ we have
$tl(g(e_{j}))>\tilde{\theta}_{e_{i}-e_{j}}(g(e_{i}))$ for $j<i\leq n$.
For the additive principal number $\tilde{\theta}_{-(e_{i}-e_{j})}(tl(g(e_{j}))$ 
we obtain
$\tilde{\theta}_{-(e_{i}-e_{j})}(tl(g(e_{j})))\geq g(e_{i})+1$
by Proposition \ref{prp:tht-}.\ref{prp:zigzag1}, and 
Proposition \ref{prp:tht-}.\ref{prp:tht6} yields
$tl(g(e_{j}))\geq\tilde{\theta}_{e_{i}-e_{i}}(g(e_{i})+1)$.
Define recursively ordinals $\xi_{i}$
by
$\xi_{n}=g(e_{n})$ and 
$\xi_{i}=g(e_{i})+\tilde{\theta}_{e_{i+1}-e_{i}}(\xi_{i+1}+1)$ for $i<n$.
We show by induction on $n-i$ that
\begin{equation}\label{eq:oflx.05.1}
\tilde{\theta}_{-(e_{i}-e_{i-1})}(tl(g(e_{i-1})))> \xi_{i}
\end{equation}
The case $i=n$ follows from $g(e_{n})=\xi_{n}$.
Assume $\tilde{\theta}_{-(e_{i+1}-e_{i})}(tl(g(e_{i})))> \xi_{i+1}$.
Since $tl(g(e_{i}))>\tilde{\theta}_{e_{i+1}-e_{i}}(g(e_{i+1}))$,
$\tilde{\theta}_{-(e_{i+1}-e_{i})}(tl(g(e_{i})))$ is additive principal, and hence
$\tilde{\theta}_{-(e_{i+1}-e_{i})}(tl(g(e_{i})))> \xi_{i+1}+1$.
Proposition \ref{prp:tht-}.\ref{prp:tht6} yields
$tl(g(e_{i}))>\tilde{\theta}_{e_{i+1}-e_{i}}( \xi_{i+1}+1)$.
Therefore
$\tilde{\theta}_{-(e_{i}-e_{i-1})}(tl(g(e_{i-1})))> 
\xi_{i}=g(e_{i})+\tilde{\theta}_{e_{i+1}-e_{i}}(\xi_{i+1}+1)$ follows from
$\tilde{\theta}_{-(e_{i}-e_{i-1})}(tl(g(e_{i-1})))> g(e_{i})$ for the additive principal number
$\tilde{\theta}_{-(e_{i}-e_{i-1})}(tl(g(e_{i-1})))$.
(\ref{eq:oflx.05.1}) is shown.

(\ref{eq:oflx.05.1}) yields
$\tilde{\theta}_{-(e-d)}(tl(g(d)))> \xi_{1}$ for $e=e_{1}$ and $e_{0}=d$.
Proposition \ref{prp:ofb} yields
$\tilde{\theta}_{-(e-d)}(\omega\cdot tl(a(g(d)))> a(\xi_{1})$, and
$\tilde{\theta}_{-(e-d)}(\omega\cdot tl(a(g(d)))> \omega\cdot a(\xi_{1})+\omega$.
We obtain $\omega\cdot tl(a(g(d))> 
\tilde{\theta}_{e-d}(\omega\cdot a(\xi_{1})+\omega)$
by Proposition \ref{prp:tht-}.\ref{prp:tht6}.

On the other side we see inductively from $a(1)=1$ and Proposition \ref{prp:ofb}
 that $o(g;e_{i})\leq \omega\cdot a(\xi_{i})$.
Therefore $\omega\cdot tl(a(g(d))> \tilde{\theta}_{e-d}(o(g;e)+\omega)$.
\\

\noindent
\ref{prp:oflx}.\ref{prp:oflx.0} by induction on the cardinality of the set 
$(d,c)\cap{\rm supp}(g)$.
Let $e=\min\{e>d: e\in{\rm supp}(g)\}$.
We obtain 
$\tilde{\theta}_{e-d}(o(g;e)+\omega)<\omega\cdot a(g(d))+\tilde{\theta}_{e-d}(o(g;e)+1)=o(g;d)$
by Proposition \ref{prp:oflx}.\ref{prp:oflx.05}.
We are done when $e\geq c$. Assume $e<c$.
IH yields $\tilde{\theta}_{c-e}(o(g;c)+1)<o(g;e)$, and hence
$\tilde{\theta}_{c-d}(o(g;c)+1)<\tilde{\theta}_{e-d}(o(g;e))<o(g;d)$.
\\

\noindent
\ref{prp:oflx}.\ref{prp:oflx.1}.
Let 
$\{d=d_{0}<d_{1}<\cdots<d_{n}\}=\{d\}\cup \supp(f^{d})$.
$f<^{d}g(d)$ means that there exists a sequence $(\xi_{i})_{i\leq n}$ of ordinals such that
$\xi_{0}$ is the shortest part of $g(d_{0})$ such that $f(d_{0})<\xi_{0}$,
and each $\xi_{i+1}$ is the shortest part of $\tilde{\theta}_{-(d_{i+1}-d_{i})}(tl(\xi_{i}))$ 
such that $f(d_{i+1})<\xi_{i+1}$ for $i<n$.
We show by induction on $n-i$ that
\begin{equation}\label{eq:of1}
d_{i}\in \supp(f) \Rarw o(f;d_{i})+\omega\leq\omega\cdot a(\xi_{i})
\end{equation}
First $o(f;d_{n})+\omega=\omega\cdot a(f(d_{n}))+\omega\leq\omega\cdot a(\xi_{n})$ by $f(d_{n})+1\leq \xi_{n}$ and
Proposition \ref{prp:ofb}.\ref{prp:ofb1}.
Next let $i<n$ with $d_{i}\in \supp(f)$. We have 
$\omega\cdot a(f(d_{i}))<\omega\cdot a(\xi_{i})$.
On the other hand we have 
$o(f;d_{i+1})+\omega\leq\omega\cdot a(\xi_{i+1})\leq 
\omega\cdot a(\tilde{\theta}_{-(d_{i+1}-d_{i})}(tl(\xi_{i})))\leq
\tilde{\theta}_{-(d_{i+1}-d_{i})}(tl(\omega\cdot a(\xi_{i})))$
by IH and Propositions \ref{prp:ofb}.\ref{prp:ofb2} and \ref{prp:ofb}.\ref{prp:ofb3}.
Proposition \ref{prp:tht-}.\ref{prp:tht6} yields
$\tilde{\theta}_{d_{i+1}-d_{i}}(o(f;d_{i+1})+1)<
\tilde{\theta}_{d_{i+1}-d_{i}}(\tilde{\theta}_{-(d_{i+1}-d_{i})}(tl(\omega\cdot a(\xi_{i}))))
\leq tl(\omega\cdot a(\xi_{i}))$.
Therefore
$o(f;d_{i})+\omega=
\omega\cdot a(f(d_{i}))+\tilde{\theta}_{d_{i+1}-d_{i}}(o(f;d_{i+1})+1)+\omega<
\omega\cdot  a(\xi_{i})$
for the additive principal number $tl(\omega\cdot a(\xi_{i}))$.
(\ref{eq:of1}) is shown.

Next we show
$o(f;d)+\omega\leq o(g;d)$.
(\ref{eq:of1}) yields 
$o(f;d)+\omega\leq\omega\cdot a(\xi_{0})\leq \omega\cdot a(g(d))\leq o(g;d)$ if $d\in \supp(f)$.
Otherwise we have 
$o(f;d_{1})+1<\omega\cdot a(\xi_{1})\leq
\tilde{\theta}_{-(d_{1}-d)}(tl(\omega\cdot a(\xi_{0})))
\leq \tilde{\theta}_{-(d_{1}-d)}(\omega\cdot a(g(d)))$ by Proposition \ref{prp:tht-}.\ref{prp:tht4}.
It follows that
$o(f;d)+\omega=\tilde{\theta}_{d_{1}-d}(o(f;d_{1})+1)+\omega<
\tilde{\theta}_{d_{1}-d}(\tilde{\theta}_{-(d_{1}-d)}(\omega\cdot a(g(d))))\leq \omega\cdot a(g(d))\leq o(g;d)$.
\\

\noindent
\ref{prp:oflx}.\ref{prp:oflx.2}.
Note that $f(c)<g(c)+\tilde{\theta}_{d-c}(g(d))\cdot \omega$ holds if 
$tl(\nu)\leq\tilde{\theta}_{d-c}(g(d))$ for a part $\nu$ of $f(c)$ with $\nu>g(c)$.
From this we see that 
$\omega\cdot a(f(c))+\omega\leq
\omega\cdot a(g(c))+\tilde{\theta}_{d-c}(\omega\cdot a(g(d)))\cdot \omega$.
We obtain
$\omega\cdot a(g(c))+\tilde{\theta}_{d-c}(\omega\cdot a(g(d)))\cdot \omega
\leq o(g;c)$ by $\omega\cdot a(g(d))\leq o(g;d)$.
Hence we can assume ${\rm supp}(f^{c+1})\neq\emptyset$.
We have $o(f;c)=\omega\cdot a(f(c))+\tilde{\theta}_{e-c}(o(f;e)+1)$
for $e=\min\{e>c: e\in {\rm supp}(f)\}$.
If $e\geq d$, then 
$\tilde{\theta}_{e-c}(o(f;e)+1)<
\tilde{\theta}_{d-c}(o(g;d)+1)$ by $o(f;d)<o(g;d)$.
Hence
$o(f;c)+\omega<\omega\cdot a(g(c))+\tilde{\theta}_{d-c}(o(g;d)+1)$.
On the other hand we have
$\tilde{\theta}_{d-c}(o(g;d)+1)\leq o(g;c)$ by Proposition \ref{prp:oflx}.\ref{prp:oflx.0}.
Therefore $o(f;c)+\omega<o(g;c)$.

In what follows assume $e<d$. Then
$f(e)<\tilde{\theta}_{d-e}(g(d))$ by 
Proposition \ref{prp:psinucomparison}.\ref{prp:psinucomparison.1}.
Next we show by induction on the cardinality of the set
$\{x\in{\rm supp}(f):  e<x<d\}$ 
that 
\begin{equation}\label{eq:oflx.2}
e\in{\rm supp}(f)\cap (c,d] \Rightarrow o(f;e)<\tilde{\theta}_{d-e}(o(g;d)+1)
\end{equation}
Let $e\in{\rm supp}(f)$ with $c<e<d$. 
If $e=\max({\rm supp}(f))$, then 
$o(f;e)=\omega\cdot a(f(e))<\tilde{\theta}_{d-e}(\omega\cdot a(g(d)))$
with $\omega\cdot a(g(d))\leq o(g;d)$.
Otherwise let $x=\min\{x\in{\rm supp}(f): x>e\}$.
Then
$o(f;e)=\omega\cdot a(f(e))+\tilde{\theta}_{x-e}(o(f;x)+1)$.
If $x\geq d$, then 
$\tilde{\theta}_{x-e}(o(f;x)+1)<
\tilde{\theta}_{d-e}(o(g;d)+1)$ by $o(f;d)<o(g;d)$.
Let $x<d$. IH yields
$o(f;x)<\tilde{\theta}_{d-x}(o(g;d)+1)$, and
$\tilde{\theta}_{x-e}(o(f;x)+1)<
\tilde{\theta}_{d-e}(o(g;d)+1)$ for the additive principal number 
$\tilde{\theta}_{d-x}(o(g;d)+1)$.
 (\ref{eq:oflx.2}) is shown.
 
 (\ref{eq:oflx.2}) yields
$\omega\cdot a(f(c)),\tilde{\theta}_{e-c}(o(f;e)+1),\omega<\tilde{\theta}_{d-c}(o(g;d)+1)$,
and $o(f;c)+\omega<\tilde{\theta}_{d-c}(o(g;d)+1)\leq o(g;c)$.
\\

\noindent
\ref{prp:oflx}.\ref{prp:oflx.3}.
We obtain $\tilde{\theta}_{c-d}(o(f;c)+1)+\omega<\tilde{\theta}_{c-d}(o(g;c))<o(g;d)$ by
Proposition \ref{prp:oflx}.\ref{prp:oflx.0}.
If $d\not\in{\rm supp}(f)$, then $o(f;d)=\tilde{\theta}_{c-d}(o(f;c)+1)$, and we are done.
Let $d\in{\rm supp}(f)$. Then
$o(f;d)=\omega\cdot a(f(d))+\tilde{\theta}_{c-d}(o(f;c)+1)$.
We obtain $\omega\cdot a(f(d))<\omega\cdot a(g(d))$.
Let $e=\min\{e>d: e\in{\rm supp}(g)\}$.
Then $o(g;d)=\omega\cdot a(g(d))+\tilde{\theta}_{e-d}(o(g;e)+1)$.
It suffices to show $\tilde{\theta}_{e-d}(o(g;e)+1)\geq\tilde{\theta}_{c-d}(o(g;c))$,
which follows from Proposition \ref{prp:oflx}.\ref{prp:oflx.0} when $e<c$.
\\

\noindent
\ref{prp:oflx}.\ref{prp:oflx.4}.
It suffices to show $o(f;d)<o(g;d)$ for $d=\max(\{0\}\cup{\rm supp}(f_{c}))$.
We can assume $c\not\in{\rm supp}(f)\cap{\rm supp}(g)$.
If $c\not\in{\rm supp}(g)$, then
$o(f;d)=\omega\cdot a(f(d))+\tilde{\theta}_{c-d}(o(f;c)+1)$ and
$o(g;d)=\omega\cdot a(f(d))+\tilde{\theta}_{e-d}(o(g;e)+1)$
for $e=\min({\rm supp}(g^{c}))$.
By $\tilde{\theta}_{e-d}(o(g;e)+1)=\tilde{\theta}_{c-d}(o(g;c))$ and
$o(f;c)+1<o(g;c)$ we obtain $o(f;d)<o(g;d)$.
Otherwise $o(f;d)=\omega\cdot a(f(d))+\tilde{\theta}_{c-d}(o(f;c))$
and $o(g;d)=\omega\cdot a(f(d))+\tilde{\theta}_{c-d}(o(g;c)+1)$.
Hence $o(f;d)<o(g;d)$.
\hspace*{\fill} $\Box$

\begin{lemma}\label{lem:of}
Let $f$ be an irreducible finite function defined from an irreducible function $g$ and ordinals $c,d$
as follows.
$f_{c}=g_{c}$, 
$c<d\in \supp(g)$ with
$(c,d)\cap \supp(g)=
(c,d)\cap \supp(f)=
\emptyset$, 
$f(c)<g(c)+\tilde{\theta}_{d-c}(g(d))\cdot\omega$, and
$f<^{d}g(d)$, cf.\,Definition \ref{df:notationsystem}.\ref{df:notationsystem.11}.
Then $o(f)<o(g)$ holds.
\end{lemma}
{\it Proof}.\hspace{2mm}
By Proposition \ref{prp:oflx}.\ref{prp:oflx.1} with $f<^{d}g(d)$ we obtain
$o(f;d)+\omega\leq o(g;d)$.
Then $o(f;c)+\omega\leq o(g;c)$ by Proposition \ref{prp:oflx}.\ref{prp:oflx.2}.
It follows that $o(f)<o(g)$ by Proposition \ref{prp:oflx}.\ref{prp:oflx.4}.
\hspace*{\fill} $\Box$

\begin{lemma}\label{lem:oflx}
For irreducible finite functions $f$ and $g$, assume
$f<_{lx}^{0}g$.
Then $o(f)<o(g)$ holds.
\end{lemma}
{\it Proof}.\hspace{2mm}
Let $b=\min\{b: f(b)\neq g(b)\}$.
It suffices to show $o(f;b)+\omega\leq o(g;b)$ by Proposition \ref{prp:oflx}.\ref{prp:oflx.4}.
According to Definition \ref{df:lx}.\ref{df:lxx}, $f<_{lx}^{0}g$ holds for $f\neq g$ iff
there exists an ordinal $d_{n}\in{\rm supp}(g)$ such that 
$f<^{d_{n}}g(d_{n})$, and one of the following holds.
\\
{\bf Case 1}.
$f(b)>g(b)$: Then there are ordinals 
$b=c_{0}<d_{1}<c_{1}<\cdots<d_{n-1}<c_{n-1}<d_{n}\,(n\geq 1)$ such that
$\{c_{i}:i<n\}\subset{\rm supp}(f)$, $\{d_{i}: 0<i<n\}\subset{\rm supp}(g)$ and
for each $c_{i}$ and $d_{i}$ the following holds.
\begin{enumerate}
\item[(a)]
$f(c_{i})>g(c_{i})$ and let $\nu_{i}$ be the shortest part of $f(c_{i})$ 
such that $\nu_{i}>g(c_{i})$. 
Then $tl(\nu_{i})\leq\tilde{\theta}_{d_{i+1}-c_{i}}(g(d_{i+1}))$
for $0\leq i<n$.

\item[(b)]
$f(d_{i})<g(d_{i})$ and let $\mu_{i}$ be the shortest part of $g(d_{i})$ 
such that $f(d_{i})<\mu_{i}$. 
Then $tl(\mu_{i})\leq\tilde{\theta}_{c_{i}-d_{i}}(f(c_{i}))$
for $0<i<n$.

\end{enumerate}
{\bf Case 2}. $f(b)<g(b)$: There are ordinals 
$b=d_{1}<c_{1}<\cdots<d_{n-1}<c_{n-1}<d_{n}\,(n\geq 1)$ such that
$\{c_{i}:0<i<n\}\subset{\rm supp}(f)$, $\{d_{i}: 0<i<n\}\subset{\rm supp}(g)$ and
for each $c_{i}\,(0<i<n)$ and $d_{i}\,(0<i<n)$ the same (a) and (b) hold.

First by Proposition \ref{prp:oflx}.\ref{prp:oflx.1} with $f<^{d_{n}}g(d_{n})$ we obtain
$o(f;d_{n})+\omega\leq o(g;d_{n})$.
If {\bf Case 2} occurs with $n=1$, then we are done.
Otherwise 
$o(f;c_{n-1})+\omega\leq o(g;c_{n-1})$ by Proposition \ref{prp:oflx}.\ref{prp:oflx.2}.
Assuming $d_{n-1}$ exists, we obtain
$o(f;d_{n-1})+\omega\leq o(g;d_{n-1})$ by Proposition \ref{prp:oflx}.\ref{prp:oflx.3}.
We see inductively $o(f;b)+\omega\leq o(g;b)$
from Propositions \ref{prp:oflx}.\ref{prp:oflx.2} and 
\ref{prp:oflx}.\ref{prp:oflx.3}.
\hspace*{\fill} $\Box$

\subsection{Computable notation system $OT_{N}$}\label{subsec:decidable}

In this subsection (except Proposition \ref{prp:l.5.4})
we work in a weak fragment of arithmetic, e.g., in the fragment $I\Sigma_{1}$ or even in the bounded arithmetic $S^{1}_{2}$.

Referring Proposition \ref{prp:psicomparison},
a set $OT_{N}$ of ordinal terms
over $\{0,\mathbb{S},+,\varphi, \Omega, \psi\}$
is defined recursively.

Simultaneously we define a subset $\Psi_{N}\subset OT_{N}$,
finite sets $K_{X}(\alpha)\subset OT_{N}$ for $X\subset OT_{N}$ and 
$\alpha\in OT_{N}$, 
cf.\,Proposition \ref{prp:l.5.4}.
$K_{\delta}(\alpha):=K_{X}(\alpha)$ with $X=\{\beta\in OT_{N}: \beta<\delta\}$.
For $\alpha=\psi_{\pi}^{f}(a)$ and $c\in \supp(f)$, let $m_{c}(\alpha)=f(c)$.
$m_{c}(\alpha)=0$ for $c\not\in \supp(f)$, and
$m(\alpha)=f$.
$E_{\mathbb{S}}(\alpha)$ denotes a finite set of subterms$<\mathbb{S}$ of $\alpha$.

For $\{\alpha_{0},\ldots,\alpha_{m},\beta\}\subset OT_{N}$,
$K_{X}(\alpha_{0},\ldots,\alpha_{m}):=\bigcup_{i\leq m}K_{X}(\alpha_{i})$,
\\
$K_{X}(\alpha_{0},\ldots,\alpha_{m})<\beta:\Leftrightarrow
\forall\gamma\in K_{X}(\alpha_{0},\ldots,\alpha_{m})(\gamma<\beta)$,
and 
\\
$\beta\leq K_{X}(\alpha_{0},\ldots,\alpha_{m}):\Leftrightarrow \exists\gamma\in K_{X}(\alpha_{0},\ldots,\alpha_{m})(\beta\leq\gamma)$.

An ordinal term in $OT_{N}$ is said to be a \textit{regular} term if it is one of the form $\mathbb{S}$,
$\Omega_{\beta}$ or $\psi_{\pi}^{f}(a)$ with the non-empty function $f$.

$\alpha=_{NF}\alpha_{m}+\cdots+\alpha_{0}$ means that $\alpha=\alpha_{m}+\cdots+\alpha_{0}$ 
and $\alpha_{m}\geq\cdots\geq\alpha_{0}$
and each $\alpha_{i}$ is a non-zero additive principal number.
$\alpha=_{NF}\varphi\beta\gamma$ means that $\alpha=\varphi\beta\gamma$ and $\beta,\gamma<\alpha$.
The ordinal function $\tilde{\theta}$ is defined in Definition \ref{df:Lam} through $\varphi$.

\begin{definition}\label{df:notationsystem}
{\rm
$\mathbb{K}:=\Omega_{\mathbb{S}+N}$ for a positive integer $N$.
$\ell\alpha$ denotes the number of occurrences of symbols
$\{0,\mathbb{S},+,\varphi, \Omega, \psi\}$
in terms $\alpha\in OT_{N}$.

 \benu
 \item\label{df:notationsystem.2}
$\{0,\Omega_{1},\mathbb{S}\}\cup\{\Omega_{\mathbb{S}+n}: 0<n\leq N\}\subset OT_{N}$.
$m(\alpha)=K_{X}(\alpha)=E_{\mathbb{S}}(\alpha)=\emptyset$ for
$\alpha\in\{0,\Omega_{1},\mathbb{S}\}\cup\{\Omega_{\mathbb{S}+n}: 0<n\leq N\}$.

 \item\label{df:notationsystem.4}
If $\alpha=_{NF}\alpha_{m}+\cdots+\alpha_{0}\, (m>0)$ with $\{\alpha_{i}:i\leq m\}\subset OT_{N}$, 
then
$\alpha\in OT_{N}$, and $m(\alpha)=\emptyset$.
$K_{X}(\alpha)=K_{X}(\alpha_{0},\ldots,\alpha_{m})$.
$E_{\mathbb{S}}(\alpha)=E_{\mathbb{S}}(\alpha_{0},\ldots,\alpha_{m})$.

 \item\label{df:notationsystem.6}
If $\alpha=_{NF}\varphi\beta\gamma$ with $\{\beta,\gamma\}\subset OT_{N}$, then
$\alpha\in OT_{N}$, and $m(\alpha)=\emptyset$.
$K_{X}(\alpha)=K_{X}(\beta,\gamma)$.
$E_{\mathbb{S}}(\alpha)=E_{\mathbb{S}}(\beta,\gamma)$.

 \item\label{df:notationsystem.8}
If $\kappa\in\Psi_{N}$, then
$\Omega_{\kappa+n}\in OT_{N}$ for $0<n\leq N$.
$m(\Omega_{\kappa+n})=\emptyset$.
$K_{X}(\Omega_{\kappa+n})=\emptyset$ if $\Omega_{\kappa+n}\in X$.
$K_{X}(\Omega_{\kappa+n})=\Gamma_{\mathbb{K}+1}$ otherwise.
$E_{\mathbb{S}}(\Omega_{\kappa+n})=\{\Omega_{\kappa+n}\}$.

 \item\label{df:notationsystem.9}
Let $\{\pi,a\}\subset OT_{N}$, $\pi$ a regular term, $K_{\alpha}(\pi,a)<a$,
and if $\pi=\Omega_{\kappa+n}$ with $\kappa\in\Psi_{N}$ and $0<k\leq N$, then 
$a<\Gamma_{\Omega_{\kappa+N}+1}$.
Then
$\alpha=\psi_{\pi}(a)\in OT_{N}$.

Let
$m(\alpha)=\emptyset$.
$K_{X}(\psi_{\pi}(a))=\emptyset$ if $\alpha\in X$.
$K_{X}(\psi_{\pi}(a))=\{a\}\cup K_{X}(a,\pi)$ if $\alpha\not\in X$ and
$\pi\in\{\Omega_{1},\mathbb{S}\}\cup\{\Omega_{\mathbb{S}+n} :0<n\leq N\}$.
$K_{X}(\psi_{\pi}(a))=\Gamma_{\mathbb{K}+1}$ otherwise.
$E_{\mathbb{S}}(\psi_{\pi}(a))=\{\psi_{\pi}(a)\}$ if $\pi\leq\mathbb{S}$.
Otherwise $E_{\mathbb{S}}(\psi_{\pi}(a))=E_{\mathbb{S}}(a)$.

 \item\label{df:notationsystem.10}
Let $\xi,a,c\in OT_{N}$,
$\xi>0$, $c<\mathbb{K}$ and
$K_{\alpha}(\xi,a,c)<a$.
Then
$\alpha=\psi_{\mathbb{S}}^{f}(a)\in \Psi_{N}$ for $f=m(\alpha)$
with ${\rm supp}(f)=\{c\}$ and $f(c)=\xi$.

Let
$K_{X}(\psi_{\mathbb{S}}^{f}(a))=\emptyset$ if $\alpha\in X$.
$K_{X}(\psi_{\mathbb{S}}^{f}(a))=\{a\}\cup K_{X}(a, c,\xi)$ otherwise.
$E_{\mathbb{S}}(\alpha)=\{\alpha\}$.

 \item\label{df:notationsystem.11}
Let $\{a,d\}\subset OT_{N}$, $\pi\in\Psi_{N}$,
$f=m(\pi)$,
 $d<c\in \supp(f)$,
and $(d,c)\cap \supp(f)=\emptyset$.

Let $g$ be an irreducible function such that 
$K(g)=\bigcup\{\{c,g(c)\}: c\in {\rm supp}(g)\}\subset OT_{N}$,
$g_{d}=f_{d}$, $(d,c)\cap \supp(g)=\emptyset$
$g(d)<f(d)+\tilde{\theta}_{c-d}(f(c))\cdot\omega$, 
and $g<^{c}f(c)$.

Then 
$\alpha=\psi_{\pi}^{g}(a)\in \Psi_{N}$ if 
$K_{\alpha}(\pi,a)\cup K_{\alpha}(K(f)\cup K(g))<a$, and, cf.\,Proposition \ref{prp:psimS},

\begin{equation}\label{eq:notationsystem.11}
E_{\mathbb{S}}(K(g))<\alpha
\end{equation}
$K_{X}(\psi_{\pi}^{g}(a))=\emptyset$ if $\alpha\in X$.
Otherwise
$K_{X}(\psi_{\pi}^{g}(a))=
\{a\}\cup K_{X}(a,\pi)\cup\bigcup\{K_{X}(b): b\in K(g)\}$.
$E_{\mathbb{S}}(\alpha)=\{\alpha\}$.

\eenu
}
\end{definition}

\begin{proposition}\label{prp:l.5.4}
For any $\alpha\in OT_{N}$ and any $X\subset OT_{N}$,
$\alpha\in\mathcal{H}_{\gamma}(X) \Lrarw K_{X}(\alpha)<\gamma$.
\end{proposition}
{\it Proof}.\hspace{2mm}
By induction on $\ell\alpha$.
\hspace*{\fill} $\Box$

\begin{lemma}\label{lem:compT}
$(OT_{N},<)$ is a computable notation system of ordinals.
Specifically
each of $\alpha<\beta$ and $\alpha=\beta$ is decidable for $\alpha,\beta\in OT_{N}$, 
and $\alpha\in OT_{N}$ is decidable for terms $\alpha$
over symbols $\{0,\mathbb{S},+,\varphi, \Omega, \psi\}$.

In particular the order type of the initial segment $\{\alpha\in OT_{N}: \alpha<\Omega_{1}\}$
is less than $\omega_{1}^{CK}$. 
\end{lemma}
{\it Proof}.\hspace{2mm}
This is seen from Propositions \ref{prp:psicomparison} and \ref{prp:l.5.4}.
\hspace*{\fill} $\Box$
\\

In what follows by ordinals we mean ordinal terms in $OT_{N}$ for a fixed positive integer $N$.

\section{Well-foundedness proof with the maximal distinguished set}\label{sec:distinguished}
In this section working in the second order arithmetic 
$\Sigma^{1-}_{2}\mbox{{\rm -CA}}+\Pi^{1}_{1}\mbox{{\rm -CA}}_{0}$, 
we show the well-foundedness of the notation system $OT_{N}$ for $N$.
The proof is based on distinguished classes, which was first
introduced by Buchholz\cite{Buchholz75}.

\subsection{Coefficients}
In this subsection we introduce coefficient sets $\mathcal{E}(\alpha),G_{\kappa}(\alpha), F_{\delta}(\alpha),k_{\delta}(\alpha)$ of 
$\alpha\in OT_{N}$,
each of which is a finite set of subterms of $\alpha$.
These are utilized in our well-foundedness proof.
Roughly $\mathcal{E}(\alpha)$ is the set of subterms of the form $\psi_{\pi}^{f}(a)$, and
$F_{\delta}(\alpha)$ [$k_{\delta}(\alpha)$] the set of subterms$<\delta$ [subterms$\geq\delta$], resp.
$G_{\kappa}(\alpha)$ is an analogue of sets $K_{\kappa}\alpha$ in \cite{odMahlo}.

\begin{definition}\label{df:prec}
{\rm
Let $pd(\psi_{\pi}^{f}(a))=\pi$ (even if $f=\emptyset$).
Moreover for $n$,
$pd^{(n)}(\alpha)$ is defined recursively by $pd^{(0)}(\alpha)=\alpha$ and
$pd^{(n+1)}(\alpha)\simeq pd(pd^{(n)}(\alpha))$.
$pd(\mathbb{S})$ is undefined.

For terms $\pi,\kappa\in OT_{N}$,
$\pi\prec\kappa$ denotes the transitive closure of the relation
$\{(\pi,\kappa): \exists f\exists b[\pi=\psi_{\kappa}^{f}(b)]\}$,
 and its reflexive closure
$\pi\preceq\kappa:\Lrarw \pi\prec\kappa \lor \pi=\kappa$.
}
\end{definition}

\begin{definition}\label{df:EGFk}

{\rm 
For terms $\alpha,\kappa,\delta\in OT_{N}$, finite sets 
$\mathcal{E}(\alpha),  G_{\kappa}(\alpha), F_{\delta}(\alpha), k_{\delta}(\alpha)\subset OT_{N}$ 
 are defined recursively as follows.
\benu

\item
$\mathcal{E}(\alpha)=\emptyset$ for $\alpha\in\{0,\Omega_{1},\mathbb{S}\}\cup\{\Omega_{\mathbb{S}+n}:0<n\leq N\}$.
$\mathcal{E}(\alpha_{m}+\cdots+\alpha_{0})=\bigcup_{i\leq m}\mathcal{E}(\alpha_{i})$.
$\mathcal{E}(\varphi\beta\gamma)=\mathcal{E}(\beta)\cup\mathcal{E}(\gamma)$.
$\mathcal{E}(\Omega_{\alpha+n})=\{\alpha\}$ for $\alpha\in\Psi_{N}$.
$\mathcal{E}(\psi_{\pi}^{f}(a))=\{\psi_{\pi}^{f}(a)\}$.

\item
$\mathcal{A}(\alpha)=\bigcup\{\mathcal{A}(\beta): \beta\in\mathcal{E}(\alpha)\}$
 for $\mathcal{A}\in\{G_{\kappa},F_{\delta},k_{\delta}\}$.

\item
\[
G_{\kappa}(\psi_{\pi}^{f}(a))=\left\{
\begin{array}{ll}
G_{\kappa}(\{\pi,a\}\cup K(f)) & \kappa<\pi
\\
G_{\kappa}(\pi) & \pi<\kappa \spand \pi\not\preceq\kappa
\\
\{\psi_{\pi}^{f}(a)\} & \pi\preceq\kappa
\end{array}
\right.
\]
\[
F_{\delta}(\psi_{\pi}^{f}(a))=\left\{
\begin{array}{ll}
F_{\delta}(\{\pi,a\}\cup K(f)) & \psi_{\pi}^{f}(a)\geq\delta
\\
\{\psi_{\pi}^{f}(a)\} & \psi_{\pi}^{f}(a)<\delta
\end{array}
\right.
\]
\[
k_{\delta}(\psi_{\pi}^{f}(a))=\left\{
\begin{array}{ll}
\{\psi_{\pi}^{f}(a)\}\cup k_{\delta}(\{\pi,a\}\cup K(f)) & \psi_{\pi}^{f}(a)\geq\delta
\\
\emptyset & \psi_{\pi}^{f}(a)<\delta
\end{array}
\right.
\]
\eenu

For $\mathcal{A}\in\{K_{\delta},G_{\kappa},F_{\delta},k_{\delta}\}$ and sets $X\subset OT_{N}$,
$\mathcal{A}(X):=\bigcup\{\mathcal{A}(\alpha): \alpha\in X\}$.
}
\end{definition}

\begin{definition}
{\rm $S(\eta)$ denotes the set of immediate subterms of $\eta$ when $\eta\not\in\mathcal{E}(\eta)$.
For example $S(\varphi\beta\gamma)=\{\beta,\gamma\}$.
$S(\eta):=\emptyset$ when $\eta\in\{0,\Omega_{1},\mathbb{S}\}\cup\{\Omega_{\mathbb{S}+n}:0<n\leq N\}$,
$S(\Omega_{\eta+n})=S(\eta)=\{\eta\}$ when $\eta\in\mathcal{E}(\eta)$.
}
\end{definition}

\begin{proposition}\label{prp:G}
For $\alpha,\kappa,a,b\in OT_{N}$,
\benu
\item\label{prp:G1}
$G_{\kappa}(\alpha)\leq\alpha$.

\item\label{prp:G2}
$\alpha\in\mathcal{H}_{a}(b) \Rarw G_{\kappa}(\alpha)\subset\mathcal{H}_{a}(b)$.

\eenu
\end{proposition}
{\it Proof}.\hspace{2mm}
These are shown simultaneously by induction on $\ell\alpha$.
It is easy to see that
\begin{equation}\label{eq:G}
G_{\kappa}(\alpha)\ni\beta \Rarw \beta\prec\kappa \spand \ell\kappa<\ell\beta\leq \ell\alpha
\end{equation}
\ref{prp:G}.\ref{prp:G1}.
Consider the case $\alpha=\psi_{\pi}^{f}(a)$ with $\pi\not\preceq\kappa$.
First let $\kappa<\pi$. Then $G_{\kappa}(\alpha)=G_{\kappa}(\{\pi,a\}\cup K(f))$.
On the other hand we have  $\forall \gamma\in K(f)\cup\{\pi,a\}(K_{\alpha}(\gamma)<a)$, i.e,
$K(f)\cup\{\pi,a\}\subset\mathcal{H}_{a}(\alpha)$.
Proposition \ref{prp:G}.\ref{prp:G2} with (\ref{eq:G}) 
yields $G_{\kappa}(K(f)\cup\{\pi,a\})\subset\mathcal{H}_{a}(\alpha)\cap\kappa\subset \mathcal{H}_{a}(\alpha)\cap\pi\subset\alpha$.
Hence $G_{\kappa}(\alpha)<\alpha$.

Next let $\pi<\kappa$ and $\pi\not\preceq\kappa$. Then $G_{\kappa}(\alpha)=G_{\kappa}(\pi)$.
By IH we have $G_{\kappa}(\pi)\leq\pi$, and $G_{\kappa}(\pi)<\pi$ by $\pi\not\preceq\kappa$.
On the other hand we have  $K_{\alpha}(\pi)<a$, i.e,
$\pi\in\mathcal{H}_{a}(\alpha)$.
Proposition \ref{prp:G}.\ref{prp:G2} 
yields $G_{\kappa}(\pi)\subset\mathcal{H}_{a}(\alpha)\cap\pi\subset\alpha$.
Hence $G_{\kappa}(\alpha)<\alpha$.
\\
\ref{prp:G}.\ref{prp:G2}.
Since $G_{\kappa}(\alpha)\leq\alpha$ by Proposition \ref{prp:G}.\ref{prp:G1}, we can assume $\alpha\geq b$.
Again consider the case $\alpha=\psi_{\pi}^{f}(a)$ with $\pi\not\preceq\kappa$.
Then $K(f)\cup\{\pi,a\}\subset\mathcal{H}_{a}(b)$ and $G_{\kappa}(\alpha)\subset G_{\kappa}(K(f)\cup\{\pi,a\})$.
IH yields the lemma.
\hspace*{\fill} $\Box$

\begin{proposition}\label{prp:G4}
Let $\gamma\preceq\tau$ and $\gamma\not\prec\kappa$. Then $G_{\kappa}(\tau)\subset G_{\kappa}(\gamma)$.
\end{proposition}
{\it Proof}.\hspace{2mm}
Let  $\gamma\not\prec\kappa$.
We show $\gamma\preceq\tau \Rarw G_{\kappa}(\tau)\subset G_{\kappa}(\gamma)$
by induction on $\ell\gamma-\ell\tau$.
Let $\gamma\preceq\tau=\psi_{\pi}^{f}(a)$.
By IH we have $G_{\kappa}(\tau)\subset G_{\kappa}(\gamma)$.
On the other hand we have $G_{\kappa}(\pi)\subset G_{\kappa}(\tau)$
since $\pi\not\prec\kappa$ and $\pi=\kappa \Rarw G_{\kappa}(\pi)=\emptyset$, cf.\,(\ref{eq:G}).
\hspace*{\fill} $\Box$

\begin{proposition}\label{prp:G3}
If
$\beta\not\in\mathcal{H}_{a}(\alpha)$ and $K_{\delta}(\beta)<a$, then there exists a $\gamma\in F_{\delta}(\beta)$
such that $\mathcal{H}_{a}(\alpha)\not\ni\gamma<\delta$.
\end{proposition}
{\it Proof}.\hspace{2mm}
By induction on $\ell\beta$. Assume $\beta\not\in\mathcal{H}_{a}(\alpha)$ and $K_{\delta}(\beta)<a$.
By IH we can assume that $\beta=\psi_{\kappa}^{f}(b)$.
If $\beta<\delta$, then $\beta\in F_{\delta}(\beta)$, and $\gamma=\beta$ is a desired one.
Assume $\beta\geq\delta$. Then we obtain $K_{\delta}(\beta)=\{b\}\cup K_{\delta}(\{b,\kappa\}\cup K(f))<a$.
In particular $b<a$, and hence $\{b,\kappa\}\cup K(f)\not\subset\mathcal{H}_{a}(\alpha)$.
By IH there exists a $\gamma\in F_{\delta}(\{b,\kappa\}\cup K(f))=F_{\delta}(\beta)$
such that $\mathcal{H}_{a}(\alpha)\not\ni\gamma<\delta$.
\hspace*{\fill} $\Box$

\subsection{Distinguished sets}

In this subsection
we establish elementary facts on
distinguished classes.
Since many properties of distinguished classes are seen as in \cite{Wienpi3d, WFnonmon2}, we will omit their proofs.
$N$ denotes a fixed positive integer.
$\mathbb{K}=\Omega_{\mathbb{S}+N}$.

 $X,Y,\ldots$ range over subsets of $OT_{N}$. 
We define sets $\mathcal{C}^{\alpha}(X)\subset OT_{N}$ for $\alpha\in OT_{N}$ and
$X\subset OT_{N}$ as follows.

\begin{definition}\label{df:CX}
{\rm 
Let $\alpha,\beta\in OT_{N}$ and $X\subset OT_{N}$.
\benu
\item
$\mathcal{C}^{\alpha}(X)$ denotes the closure of $\{0,\Omega_{1},\mathbb{S}\}\cup\{\Omega_{\mathbb{S}+n} :0<n\leq N\}\cup(X\cap\alpha)$
under $+$, 
$\Psi_{N}\ni\sigma\mapsto \Omega_{\sigma+n}\,(0<n\leq N)$,
$(\beta,\gamma)\mapsto\varphi\beta\gamma$,
and $(\sigma,\beta,f)\mapsto \psi_{\sigma}^{f}(\beta)$
for $ \sigma>\alpha$ in  $OT_{N}$.

The last clause says that,
$\psi_{\sigma}^{f}(\beta)\in\mathcal{C}^{\alpha}(X)$ if
$\{\sigma,\beta\}\cup K(f)\subset\mathcal{C}^{\alpha}(X)$ and $\sigma>\alpha$.

\item
$\alpha^{+}=\Omega_{a+1}$ denotes the least regular term above $\alpha$ if such a term exists. 
Otherwise $\alpha^{+}:=\infty$.

\eenu
}
\end{definition}

\begin{proposition}\label{lem:CX2} {\rm (Cf. \cite{Wienpi3d}, Lemmas 3.5.3 and 3.5.4.)}
\\
Assume $\forall\gamma\in X[ \gamma\in\mathcal{C}^{\gamma}(X)]$ for a set $X\subset OT_{N}$.

\benu
\item\label{lem:CX2.3} 
$\alpha\leq\beta \Rarw \mathcal{C}^{\beta}(X)\subset \mathcal{C}^{\alpha}(X)$.

\item\label{lem:CX2.4} 
$\alpha<\beta<\alpha^{+} \Rarw \mathcal{C}^{\beta}(X)=\mathcal{C}^{\alpha}(X)$.
\eenu
\end{proposition}
{\it Proof}.\hspace{2mm}
\ref{lem:CX2}.\ref{lem:CX2.3}.
We see by induction on $\ell\gamma\,(\gamma\in OT_{N})$ that
\begin{equation}\label{eq:CX2.3}
\forall\beta\geq\alpha[\gamma\in \mathcal{C}^{\beta}(X) \Rarw \gamma\in \mathcal{C}^{\alpha}(X)\cup(X\cap\beta)]
\end{equation}
For example, if $\psi_{\pi}^{f}(\delta)\in \mathcal{C}^{\beta}(X)$ with $\pi>\beta\geq\alpha$
and $\{\pi,\delta\}\cup K(f)\subset\mathcal{C}^{\alpha}(X)\cup(X\cap\beta)$, then 
$\pi\in\mathcal{C}^{\alpha}(X)$, and 
for any $\gamma\in\{\delta\}\cup K(f)$, either $\gamma\in\mathcal{C}^{\alpha}(X)$ or
$\gamma\in X\cap\beta$. If $\gamma<\alpha$, then $\gamma\in X\cap\alpha\subset\mathcal{C}^{\alpha}(X)$.
If $\alpha\leq\gamma\in X\cap\beta$, then $\gamma\in\mathcal{C}^{\gamma}(X)$ by the assumption, and
by IH we have $\gamma\in\mathcal{C}^{\alpha}(X)\cup(X\cap\gamma)$, i.e., $\gamma\in\mathcal{C}^{\alpha}(X)$.
Therefore $\{\pi,\delta\}\cup K(f)\subset\mathcal{C}^{\alpha}(X)$, and 
$\psi_{\pi}^{f}(\delta)\in\mathcal{C}^{\alpha}(X)$.

Using (\ref{eq:CX2.3}) we see from the assumption that
$
\forall\beta\geq\alpha[ \gamma\in\mathcal{C}^{\beta}(X) \Rarw \gamma\in\mathcal{C}^{\alpha}(X)]
$.
\\

\noindent
\ref{lem:CX2}.\ref{lem:CX2.4}.
Assume $\alpha<\beta<\alpha^{+}$. Then by Proposition \ref{lem:CX2}.\ref{lem:CX2.3} we have
$\mathcal{C}^{\beta}(X)\subset\mathcal{C}^{\alpha}(X)$.
$\mathcal{C}^{\alpha}(X) \subset\mathcal{C}^{\beta}(X)$ is easily seen from $\beta<\alpha^{+}$.
\hspace*{\fill} $\Box$

\begin{definition}\label{df:wftg}
{\rm
\benu
\item 
$Prg[X,Y] :\Lrarw \forall\alpha\in X(X\cap\alpha\subset Y \to \alpha\in Y)$.

\item 
For a definable class $\mathcal{X}$, $TI[\mathcal{X}]$ denotes the schema:\\
$TI[\mathcal{X}] :\Lrarw Prg[\mathcal{X},\mathcal{Y}]\to \mathcal{X}\subset\mathcal{Y} \mbox{ {\rm holds for} any definable classes } \mathcal{Y}$.
\item
For $X\subset OT_{N}$, $W(X)$ denotes the \textit{well-founded part} of $X$. 
\item 
$Wo[X] : \Lrarw X\subset W(X)$.
\eenu
}
\end{definition}
Note that for $\alpha\in OT_{N}$,
$W(X)\cap\alpha=W(X\cap\alpha)$.

\begin{definition} \label{df:3wfdtg32}
{\rm 
For $X\subset OT_{N}$ and
$\alpha\in OT_{N}$,
\benu
\item\label{df:3wfdtg.832}
\begin{equation}\label{eq:distinguishedclass}
D[X] :\Lrarw 
\forall\alpha(\alpha\leq X\to W(\mathcal{C}^{\alpha}(X))\cap\alpha^{+}= X\cap\alpha^{+})
\end{equation}

A set $X$ is said to be a \textit{distinguished set} if $D[X]$.

\item\label{df:3wfdtg.932}
$\mathcal{W}:=\bigcup\{X :D[X]\}$.

\eenu
}
\end{definition}

Let $\alpha\in X$ for a distinguished set $X$. 
Then $W(\mathcal{C}^{\alpha}(X))\cap\alpha^{+}= X\cap\alpha^{+}$.
Hence $X$ is a well order.

\begin{proposition}\label{lem:3.11.632}
{\rm (Cf. Lemma 3.30 in \cite{WFnonmon2}.)}\\
Let $X$ be a distinguished set. Then
$\alpha\in X \Rarw \forall\beta[\alpha\in\mathcal{C}^{\beta}(X)]$.
\end{proposition}
{\it Proof}.\hspace{2mm}
Let $D[X]$ and $\alpha\in X$. 
Then $\alpha\in X\cap\alpha^{+}=W(\mathcal{C}^{\alpha}(X))\cap\alpha^{+}\subset \mathcal{C}^{\alpha}(X)$.
Hence $\forall\gamma\in X(\gamma\in\mathcal{C}^{\gamma}(X))$, and
$\alpha\in\mathcal{C}^{\beta}(X)$ for any $\beta\leq\alpha$ by Proposition \ref{lem:CX2}.\ref{lem:CX2.3}. 
Moreover for $\beta>\alpha$ we have $\alpha\in X\cap\beta\subset\mathcal{C}^{\beta}(X)$.
\hspace*{\fill} $\Box$

\begin{proposition}\label{lem:5uv.232general}{\rm (Cf. Lemma 3.28 in \cite{WFnonmon2}.)}\\
For any distinguished sets $X$ and $Y$, the following holds:
\[
X\cap\alpha=Y\cap\alpha \Rarw \forall\beta<\alpha^{+}\{\mathcal{C}^{\beta}(X)\cap\beta^{+}=\mathcal{C}^{\beta}(Y)\cap\beta^{+}\}.
\]
\end{proposition}
{\it Proof}.\hspace{2mm}
Assume that $X\cap\alpha=Y\cap\alpha$ and $\alpha\leq\beta<\alpha^{+}$.
We obtain $\mathcal{C}^{\alpha}(X)=\mathcal{C}^{\alpha}(Y)$.
On the other hand we have $\mathcal{C}^{\beta}(X)=\mathcal{C}^{\alpha}(X)$ and similarly for 
$\mathcal{C}^{\beta}(Y)$
by Proposition \ref{lem:CX2}.\ref{lem:CX2.4}.
Hence $\mathcal{C}^{\beta}(X)=\mathcal{C}^{\beta}(Y)$.
\hspace*{\fill} $\Box$

\begin{proposition}\label{lem:3wf5}
For distinguished sets $X$ and $Y$,

$\alpha\leq X \spand \alpha\leq Y \Rarw X\cap\alpha^+=Y\cap\alpha^+$.

\end{proposition}

\begin{proposition}\label{lem:3wf6}
$\mathcal{W}$ is the maximal distinguished {\rm set}, i.e.,
$D[\mathcal{W}]$ and $\exists X(X=\mathcal{W})$.
\end{proposition}
{\it Proof}.\hspace{2mm} 
$D[\mathcal{W}]$  is seen from
Proposition \ref{lem:3wf5}.
The axiom $\Sigma^{1-}_{2}\mbox{{\rm -CA}}$ yields $\exists X(X=\mathcal{W})$.
\hspace*{\fill} $\Box$

\begin{proposition}\label{lem:CX3} {\rm (Cf. \cite{Wienpi3d}, Lemma 3.6.)}
Let $\gamma<\beta$.
Assume $\alpha\in\mathcal{C}^{\gamma}(\mathcal{W})$ and $\forall\kappa\leq\beta[G_{\kappa}(\alpha)<\gamma]$.
Moreover assume 
$\forall\delta[\ell\delta\leq\ell\alpha\spand\delta\in\mathcal{C}^{\gamma}(\mathcal{W})\cap\gamma\Rarw\delta\in \mathcal{C}^{\beta}(\mathcal{W})]$.
Then $\alpha\in\mathcal{C}^{\beta}(\mathcal{W})$.
\end{proposition}
{\it Proof}.\hspace{2mm}
By induction on $\ell\alpha$. 
If $\alpha<\gamma$, then $\alpha\in \mathcal{C}^{\gamma}(\mathcal{W})\cap\gamma$.
The third assumption yields $\alpha\in \mathcal{C}^{\beta}(\mathcal{W})$. 
Assume $\alpha\geq\gamma$. 
Except the case $\alpha=\psi_{\pi}^{f}(a)$ for some $\pi,a,f$, IH yields $\alpha\in\mathcal{C}^{\beta}(\mathcal{W})$. 
Suppose $\alpha=\psi_{\pi}^{f}(a)$ for some $\{\pi,a\}\cup K(f)\subset\mathcal{C}^{\gamma}(\mathcal{W})$ and $\pi>\gamma$.
If $\pi\leq\beta$, then $\{\alpha\}=G_{\pi}(\alpha)<\gamma$ by the second assumption. Hence this is not the case, and we obtain 
$\pi>\beta$.
 Then $G_{\kappa}(\{\pi,a\}\cup K(f))=G_{\kappa}(\alpha)<\gamma$ for any $\kappa\leq\beta<\pi$. 
IH yields $\{\pi,a\}\cup K(f)\subset\mathcal{C}^{\beta}(\mathcal{W})$. We conclude $\alpha\in\mathcal{C}^{\beta}(\mathcal{W})$ from $\pi>\beta$.
\hspace*{\fill} $\Box$

\begin{proposition}\label{prp:updis}
$\alpha\leq \mathcal{W}\cap\beta^{+} \spand \alpha\in\mathcal{C}^{\beta}(\mathcal{W}) \Rarw \alpha\in \mathcal{W}$.
\end{proposition}
{\it Proof}.\hspace{2mm}
This is seen from Propositions \ref{lem:CX2}, \ref{lem:3.11.632} an \ref{lem:3wf6}.
\hspace*{\fill} $\Box$

\begin{proposition}\label{prp:maxwup1}
Let $\alpha\in \mathcal{W}$.
Then $S(\alpha)\subset \mathcal{W}$.
\end{proposition}
{\it Proof}.\hspace{2mm}
Let $\alpha\in \mathcal{W}$. Then 
$\alpha\in\mathcal{C}^{\alpha}(\mathcal{W})$.
Hence $S(\alpha)\cap\alpha\subset\mathcal{C}^{\alpha}(\mathcal{W})\cap\alpha\subset \mathcal{W}$.
\hspace*{\fill} $\Box$

\begin{proposition}\label{prp:id3wf20a-1}
$\alpha\in\mathcal{C}^{\delta}(\mathcal{W})\Rarw F_{\delta}(\alpha)\subset \mathcal{W}$.
\end{proposition}
{\it Proof} by induction on $\ell\alpha$.
If $\alpha\in \mathcal{W}\cap\delta$, then $S(\alpha)\subset \mathcal{W}$ by Proposition \ref{prp:maxwup1}, and 
$F_{\delta}(\alpha)=F_{\delta}(S(\alpha))\subset \mathcal{W}$ by IH.
Otherwise $S(\alpha)\subset\mathcal{C}^{\delta}(\mathcal{W})$, and $F_{\delta}(\alpha)=F_{\delta}(S(\alpha))\subset \mathcal{W}$ by IH.
\hspace*{\fill} $\Box$

\subsection{Sets $\mathcal{C}^{\alpha}(\mathcal{W})$ and $\mathcal{G}$}\label{subsec:C(X)}

In this subsection we establish a key fact, Lemma \ref{th:3wf16} on distinguished sets.

\begin{definition}\label{df:calg}
$\mathcal{G}:=\{\alpha\in OT_{N} :\alpha\in \mathcal{C}^{\alpha}(\mathcal{W}) \spand \mathcal{C}^{\alpha}(\mathcal{W})\cap\alpha\subset \mathcal{W}\}$.
\end{definition}

\begin{lemma}\label{lem:wf5.332}
$\mathcal{W}\subset\mathcal{G}$.
\end{lemma}
{\it Proof}.\hspace{2mm} 
Let $\gamma\in\mathcal{W}$.
We have $\gamma\in\mathcal{C}^{\gamma}(\mathcal{W})$.
Assume $\alpha\in\mathcal{C}^{\gamma}(\mathcal{W})\cap\gamma$. 
We have $\gamma\in W(\mathcal{C}^{\gamma}(\mathcal{W}))\cap\gamma^{+}=\mathcal{W}\cap\gamma^{+}$ by $\gamma\in \mathcal{W}$. 
Hence $\alpha\in W(\mathcal{C}^{\gamma}(\mathcal{W}))\cap\gamma^{+}\subset \mathcal{W}$, and
$\mathcal{C}^{\gamma}(\mathcal{W})\cap\gamma\subset\mathcal{W}$.
We obtain $\mathcal{W}\subset\mathcal{G}$.
\hspace*{\fill} $\Box$

\begin{proposition}\label{lem:CMsmpl}
Let $\alpha\in\mathcal{C}^{\beta}(\mathcal{W})$.
Assume $\mathcal{W}\cap\beta<\delta$.
Then $F_{\delta}(\alpha)\subset\mathcal{C}^{\beta}(\mathcal{W})$.
\end{proposition}
{\it Proof}.\hspace{2mm}
By induction on $\ell\alpha$.
Let $\{0,\Omega_{1},\mathbb{S}\}\cup\{\Omega_{\mathbb{S}+n}:0<n\leq N\}\not\ni\alpha\in\mathcal{C}^{\beta}(\mathcal{W})$.
First consider the case $\alpha\not\in\mathcal{E}(\alpha)$.
If $\alpha\in \mathcal{W}\cap\beta\subset\mathcal{G}$, then 
$\mathcal{E}(\alpha)\subset\mathcal{C}^{\alpha}(\mathcal{W})\cap\alpha\subset \mathcal{W}\subset\mathcal{C}^{\beta}(\mathcal{W})$ by Proposition \ref{lem:3.11.632}.
Otherwise we have $\alpha\not\in\mathcal{E}(\alpha)\subset\mathcal{C}^{\beta}(\mathcal{W})$.
In each case IH yields $F_{\delta}(\alpha)=F_{\delta}(\mathcal{E}(\alpha))\subset\mathcal{C}^{\beta}(\mathcal{W})$.

Let $\alpha=\psi_{\pi}^{f}(a)$ for some $\pi,f,a$. 
If $\alpha<\delta$, then $F_{\delta}(\alpha)=\{\alpha\}$, and there is nothing to prove.
Let $\alpha\geq\delta$. Then $F_{\delta}(\alpha)=F_{\delta}(\{\pi,a\}\cup K(f))$.
On the other side we see $\{\pi,a\}\cup K(f)\subset\mathcal{C}^{\beta}(\mathcal{W})$ from $\alpha\in\mathcal{C}^{\beta}(\mathcal{W})$ 
and the assumption.
IH yields $F_{\delta}(\alpha)\subset\mathcal{C}^{\beta}(\mathcal{W})$.
\hspace*{\fill} $\Box$

\begin{proposition}\label{lem:id7wf8}
\benu
\item\label{lem:id7wf8.3}
$\forall\tau[\alpha\in \mathcal{W}
\Rarw G_{\tau}(\alpha)\subset \mathcal{W}]$.

\item\label{lem:id7wf8.3.5}
$\forall\beta\forall\tau[\alpha\in\mathcal{C}^{\beta}(\mathcal{W}) \Rarw G_{\tau}(\alpha)\subset \mathcal{C}^{\beta}(\mathcal{W})]$.
\eenu
\end{proposition}
{\it Proof}.\hspace{2mm}
 By simultaneous induction on $\ell\alpha$.
\\
\ref{lem:id7wf8}.\ref{lem:id7wf8.3}. 
Suppose $\alpha\in \mathcal{W}$. 
Then $\alpha\in\mathcal{G}$ by Lemma \ref{lem:wf5.332}, i.e.,
$\alpha\in \mathcal{C}^{\alpha}(\mathcal{W})$, and $\mathcal{C}^{\alpha}(\mathcal{W})\cap\alpha\subset \mathcal{W}$.

Let $\alpha\not\in\mathcal{E}(\alpha)$. Then $\mathcal{E}(\alpha)\subset \mathcal{C}^{\alpha}(\mathcal{W})\cap\alpha\subset \mathcal{W}$.
 IH yields
$G_{\tau}(\alpha)=G_{\tau}(\mathcal{E}(\alpha))\subset \mathcal{W}$.
Assume $\alpha\in\mathcal{E}(\alpha)$, i.e., $\alpha=\psi_{\pi}^{f}(a)$ for some $\pi,a,f$.
Then $\{\pi,a\}\cup K(f)\subset\mathcal{C}^{\alpha}(\mathcal{W})$ by $\alpha\in\mathcal{C}^{\alpha}(\mathcal{W})$.
We can assume $\pi\not\preceq\tau$. Then $G_{\tau}(\alpha)\subset G_{\tau}(\{\pi,a\}\cup K(f))$.
By IH with Proposition \ref{prp:G}.\ref{prp:G1} we obtain
$G_{\tau}(\alpha)\subset G_{\tau}(\{\pi,a\}\cup K(f))\subset\mathcal{C}^{\alpha}(\mathcal{W})\cap\alpha\subset \mathcal{W}$.
\\

\noindent
\ref{lem:id7wf8}.\ref{lem:id7wf8.3.5}. 
Assume $\alpha\in\mathcal{C}^{\beta}(\mathcal{W})$. 
We show $G_{\tau}(\alpha)\subset\mathcal{C}^{\beta}(\mathcal{W})$. 
If $\alpha\in \mathcal{W}\cap\beta$, then by Proposition \ref{lem:id7wf8}.\ref{lem:id7wf8.3} we obtain
$G_{\tau}(\alpha)\subset \mathcal{W}\cap\beta\subset \mathcal{C}^{\beta}(\mathcal{W})$. 
Consider the case $\alpha\not\in \mathcal{W}\cap\beta$.
If $\alpha\not\in\mathcal{E}(\alpha)$, then IH yields $G_{\tau}(\alpha)=G_{\tau}(\mathcal{E}(\alpha))\subset\mathcal{C}^{\beta}(\mathcal{W})$.
Let $\alpha=\psi_{\pi}^{f}(a)$ for some $\{\pi,a\}\cup K(f)\subset\mathcal{C}^{\beta}(\mathcal{W})$ with 
$\beta<\pi\not\preceq\tau$.
IH yields $G_{\tau}(\alpha)\subset G_{\tau}(\{\pi,a\}\cup K(f))\subset\mathcal{C}^{\beta}(\mathcal{W})$.
\hspace*{\fill} $\Box$

\begin{proposition}\label{lem:KC} 
Assume $\alpha\in\mathcal{C}^{\beta}(\mathcal{W})$. Then
$\forall\sigma\leq\beta(G_{\sigma}(\alpha)\subset \mathcal{W})$.
\end{proposition}
{\it Proof}.\hspace{2mm} 
By induction on $\ell\alpha$ using Proposition \ref{lem:id7wf8}.\ref{lem:id7wf8.3}
 we see $\alpha\in\mathcal{C}^{\beta}(\mathcal{W})\spand\sigma\leq\beta \Rarw G_{\sigma}(\alpha)\subset \mathcal{W}$. 
\hspace*{\fill} $\Box$

\begin{proposition}\label{lem:3wf9}
Put $Y=W(\mathcal{C}^{\alpha}(\mathcal{W}))\cap\alpha^{+}$.
Assume that $\alpha\in\mathcal{G}$ and
\[
\forall\beta<\mathbb{K}(\mathcal{W}<\beta \spand \beta^+<\alpha^+ \Rarw 
W(\mathcal{C}^{\beta}(\mathcal{W}))\cap\beta^{+}\subset \mathcal{W})
.\]
Then $\alpha\in Y$ and $D[Y]$.
\end{proposition}
{\it Proof}.\hspace{2mm} 
As in \cite{odMahlo, WFnonmon2}  
this is seen from Lemma \ref{lem:wf5.332}.
\hspace*{\fill} $\Box$

\begin{proposition}\label{lem:3wf8}
$0\in \mathcal{W}$.
\end{proposition}

The following Lemma \ref{th:3wf16} is a key result on distinguished classes.

\begin{lemma}\label{th:3wf16}{\rm (Cf.\,Lemma 3.3.7 in \cite{WFnonmon2}.)}\\
Suppose 
\begin{equation}\label{eq:3wf16hyp.132}
\eta\in\mathcal{G}
\end{equation}
and 
\begin{equation}\label{eq:3wf16hyp.232}
\forall\gamma\prec\eta(\gamma\in\mathcal{G}
\to \gamma\in \mathcal{W})
\end{equation}
 Then
\[
\eta\in W(\mathcal{C}^{\eta}(\mathcal{W}))\cap\eta^{+} \mbox{ and } D[W(\mathcal{C}^{\eta}(\mathcal{W}))\cap\eta^{+}]
.\]
\end{lemma}
{\it Proof}.\hspace{2mm} 
By Proposition \ref{lem:3wf9} and the hypothesis (\ref{eq:3wf16hyp.132}) it suffices to show that
\[
\forall\beta<\mathbb{K}(\mathcal{W}<\beta \spand \beta^+<\eta^+ \Rarw 
W(\mathcal{C}^{\beta}(\mathcal{W}))\cap\beta^+\subset \mathcal{W})
.\]
Assume $\mathcal{W}<\beta<\mathbb{K}$ and $\beta^+<\eta^+$. 
We have to show  $W(\mathcal{C}^{\beta}(\mathcal{W}))\cap\beta^+\subset \mathcal{W}$. 
We prove this by induction on $\gamma\in W(\mathcal{C}^{\beta}(\mathcal{W}))\cap\beta^+$. 
Suppose $\gamma\in \mathcal{C}^{\beta}(\mathcal{W})\cap\beta^{+}$ and 
\[
\mbox{{\rm MIH : }} \mathcal{C}^{\beta}(\mathcal{W})\cap\gamma\subset \mathcal{W}.
\]
We show $\gamma\in \mathcal{W}$. 

First note that
$
\gamma\leq \mathcal{W} \Rarw \gamma\in \mathcal{W}
$
since if $\gamma\leq \delta$ for some $\delta\in \mathcal{W}$, then by $\mathcal{W}<\beta$ and $\gamma\in \mathcal{C}^{\beta}(\mathcal{W})$ 
we obtain $\delta<\beta$, 
$\gamma\in \mathcal{C}^{\delta}(\mathcal{W})$ and $\delta\in W(\mathcal{C}^{\delta}(\mathcal{W}))\cap\delta^{+}=\mathcal{W}\cap\delta^{+}$.
Hence $\gamma\in W(\mathcal{C}^{\delta}(\mathcal{W}))\cap\delta^{+}\subset \mathcal{W}$.
Therefore we can assume that
\begin{equation}
\label{eq:3wf9hyp.232X}
\mathcal{W}<\gamma
\end{equation}

We show first 
\begin{equation}
\label{eq:3wf9hyp.232}
\gamma\in\mathcal{G}
\end{equation}
First $\gamma\in \mathcal{C}^{\gamma}(\mathcal{W})$ by $\gamma\in \mathcal{C}^{\beta}(\mathcal{W})\cap\beta^{+}$ and Proposition \ref{lem:CX2}. 
Second we show the following claim by induction on $\ell\alpha$:

\begin{equation}\label{clm:3wf1632}
\alpha\in\mathcal{C}^{\gamma}(\mathcal{W})\cap\gamma \Rarw  \alpha\in \mathcal{W}
\end{equation}
{\it Proof} of (\ref{clm:3wf1632}). 
Assume $\alpha\in\mathcal{C}^{\gamma}(\mathcal{W})\cap\gamma$.
We have $\alpha\in\mathcal{C}^{\beta}(\mathcal{W})\cap\gamma\Rarw \alpha\in \mathcal{W}$ by MIH.
We can assume $\gamma^{+}\leq\beta$ for otherwise we have 
$\alpha\in \mathcal{C}^{\gamma}(\mathcal{W})\cap\gamma=\mathcal{C}^{\beta}(\mathcal{W})\cap\gamma\subset \mathcal{W}$ by MIH.
In what follows assume $\alpha\not\in \mathcal{W}$.

First consider the case $\alpha\not\in\mathcal{E}(\alpha)$. 
By induction hypothesis on lengths we have $\mathcal{E}(\alpha)\subset \mathcal{W}\subset\mathcal{C}^{\beta}(\mathcal{W})$, 
and hence $\alpha\in \mathcal{C}^{\beta}(\mathcal{W})\cap\gamma$.
Therefore $\alpha\in \mathcal{W}$ by MIH.
In what follows assume $\alpha=\psi_{\pi}^{f}(a)$ for some $\pi>\gamma$ such that
 $\{\pi,a\}\cup K(f)\subset\mathcal{C}^{\gamma}(\mathcal{W})$.
\\
{\bf Case 1}. $\beta<\pi$: 
Then $\forall\kappa\leq\beta[G_{\kappa}(\{\pi,a\}\cup K(f))=G_{\kappa}(\alpha)<\alpha<\gamma]$ by Proposition \ref{prp:G}.\ref{prp:G1}.
 Proposition \ref{lem:CX3} with induction hypothesis on lengths yields 
$\{\pi,a\}\cup K(f)\subset\mathcal{C}^{\beta}(\mathcal{W})$.
Hence $\alpha\in \mathcal{C}^{\beta}(\mathcal{W})\cap\gamma$ by $\pi>\beta$.
MIH yields $\alpha\in \mathcal{W}$.
\\
{\bf Case 2}. $\beta\geq\pi$: 
We have $\alpha<\gamma<\pi\leq\beta$. 
It suffices to show that $\alpha\leq \mathcal{W}$.
Then by (\ref{eq:3wf9hyp.232X}) we have $\alpha\leq\delta\in \mathcal{W}$ for some $\delta<\gamma$, and
$\alpha\in\mathcal{W}$ by Proposition \ref{prp:updis}.

Consider first the case $\gamma\not\in\mathcal{E}(\gamma)$.
By $\alpha=\psi_{\pi}^{f}(a)<\gamma<\pi$,
we can assume that 
$\gamma\not\in\{0,\Omega_{1},\mathbb{S}\}\cup\{\Omega_{\mathbb{S}+n}:0<n\leq N\}$.
Then let $\delta=\max S(\gamma)$ denote the largest immediate subterm of $\gamma$.
Then 
$\delta\in\mathcal{C}^{\gamma}(\mathcal{W})\cap\gamma$, and 
by (\ref{eq:3wf9hyp.232X}), $\mathcal{W}<\gamma\in\mathcal{C}^{\beta}(\mathcal{W})$ we have 
$\delta\in \mathcal{C}^{\beta}(\mathcal{W})\cap\gamma$.
Hence $\delta\in \mathcal{W}$ by MIH.
Also by $\alpha<\gamma$, we obtain $\alpha\leq\delta$, i.e., $\alpha\leq \mathcal{W}$, and we are done.

Next let $\gamma=\psi_{\kappa}^{f}(b)$ for some $b,f$ and $\kappa>\beta$ by (\ref{eq:3wf9hyp.232X}) and $\gamma\in \mathcal{C}^{\beta}(\mathcal{W})$.
We have $\alpha<\gamma<\pi\leq\beta<\kappa$.
$\pi\not\in\mathcal{H}_{b}(\gamma)$ since otherwise by $\pi<\kappa$ we would have $\pi<\gamma$.
Then by Proposition \ref{prp:psicomparison} we have $a\geq b$ and $K(f)\cup\{\kappa,b\}\not\subset\mathcal{H}_{a}(\alpha)$.
On the other hand we have $K_{\gamma}(K(f)\cup\{\kappa,b\})<b\leq a$.
By Proposition \ref{prp:G3} pick a $\delta\in F_{\gamma}(K(f)\cup\{\kappa,b\})$
such that $\mathcal{H}_{a}(\alpha)\not\ni\delta<\gamma$.
Also we have $K(f)\cup\{\kappa,b\}\subset\mathcal{C}^{\beta}(\mathcal{W})$,
and $\mathcal{W}<\gamma$ by (\ref{eq:3wf9hyp.232X}).
Therefore by Proposition \ref{lem:CMsmpl} with MIH
we obtain 
$\alpha\leq\delta\in\mathcal{C}^{\beta}(\mathcal{W})\cap\gamma\subset\mathcal{W}$.
We are done.
Thus (\ref{clm:3wf1632}) is shown.
\hspace*{\fill} $\Box$
\\

\noindent
Hence we obtain (\ref{eq:3wf9hyp.232}), $\gamma\in\mathcal{G}$.
We have $\gamma<\beta^{+}\leq\eta$ and $\gamma\in\mathcal{C}^{\gamma}(\mathcal{W})$.
If $\gamma\prec\eta$, then the hypothesis (\ref{eq:3wf16hyp.232}) yields $\gamma\in \mathcal{W}$.
In what follows assume $\gamma\not\prec\eta$.

If $\forall\tau\leq\eta[G_{\tau}(\gamma)<\gamma]$, then  Proposition \ref{lem:CX3} yields 
$\gamma\in\mathcal{C}^{\eta}(\mathcal{W})\cap\eta\subset \mathcal{W}$ by $\eta\in\mathcal{G}$.

Suppose $\exists\tau\leq\eta[G_{\tau}(\gamma)=\{\gamma\}]$. 
This means, by $\gamma\not\prec\eta$, that
$\gamma\prec\tau$ for a $\tau<\eta$.
Let $\tau$ denote the maximal such one.
We have $\gamma<\tau<\eta$.
From $\gamma\in\mathcal{C}^{\gamma}(\mathcal{W})$ we see $\tau\in\mathcal{C}^{\gamma}(\mathcal{W})$.
Next we show that
\begin{equation}\label{eq:3wf1632last}
\forall\kappa\leq\eta[G_{\kappa}(\tau)<\gamma]
\end{equation}
Let $\kappa\leq\eta$. 
If $\gamma\not\prec\kappa$, then $G_{\kappa}(\tau)\subset G_{\kappa}(\gamma)<\gamma$ by 
Propositions \ref{prp:G4} and \ref{prp:G}.\ref{prp:G1}.
If $\gamma\prec\kappa$, then by the maximality of $\tau$ we have $\kappa\preceq\tau$, and hence
$G_{\kappa}(\tau)=\emptyset$, cf.  (\ref{eq:G}).
(\ref{eq:3wf1632last}) is shown.

Proposition \ref{lem:CX3} yields 
$\tau\in\mathcal{C}^{\eta}(\mathcal{W})$, and
$\tau\in\mathcal{C}^{\eta}(\mathcal{W})\cap\eta\subset \mathcal{W}$ by $\eta\in\mathcal{G}$.
Therefore $\mathcal{W}<\gamma<\tau\in \mathcal{W}$.
This is not the case by (\ref{eq:3wf9hyp.232X}).
We are done.
\hspace*{\fill} $\Box$

\subsection{Well-foundedness proof concluded}

Recall that $\mathbb{K}=\Omega_{\mathbb{S}+N}$ for a positive integer $N$.

\begin{definition}
{\rm
Define recursively sets $\mathcal{C}_{n}, \mathcal{W}_{n}\, (n\leq N+1)$ from the maximal distinguished set $\mathcal{W}$ as follows.
$\mathcal{W}_{0}=\mathcal{W}\cap\mathbb{S}$, $\mathcal{C}_{n}=\mathcal{C}^{\Omega_{\mathbb{S}+n}}(\mathcal{W}_{n})$, and
$\mathcal{W}_{n+1}=W(\mathcal{C}_{n})\cap\Omega_{\mathbb{S}+n+1}$, where $\Omega_{\mathbb{S}+N+1}=\infty$.
}
\end{definition}

In $\Sigma^{1-}_{2}\mbox{{\rm -CA}}+\Pi^{1}_{1}\mbox{{\rm -CA}}_{0}$, each $\mathcal{W}_{n}$ exists as a set.

\begin{proposition}\label{prp:pi11.0}
\benu
\item\label{prp:pi11.0.2}
$\mathcal{C}_{n}\cap\Omega_{\mathbb{S}+n}=\mathcal{W}_{n}$ and $\Omega_{\mathbb{S}+n}\in\mathcal{W}_{n+1}$.
\item\label{prp:pi11.0.0}
$\mathcal{C}_{n+1}\subset\mathcal{C}_{n}$.
\item\label{prp:pi11.0.1}
$\mathcal{W}_{n}\subset\mathcal{W}_{n+1}$.
\item\label{prp:pi11.0.4}
$\alpha\in\mathcal{W}_{n}\Lrarw \mathcal{E}(\alpha)\subset\mathcal{W}_{n}$ for each $n$.
\item\label{prp:pi11.0.3}
$\mathcal{C}_{N}=W(\mathcal{C}_{N-1})=W(\mathcal{C}_{N})$.
\eenu
\end{proposition}
{\it Proof}.\hspace{2mm}
\ref{prp:pi11.0}.\ref{prp:pi11.0.2}.
$\mathcal{C}_{n}\cap\Omega_{\mathbb{S}+n}=\mathcal{C}^{\Omega_{\mathbb{S}+n}}(\mathcal{W}_{n})\cap\Omega_{\mathbb{S}+n}=\mathcal{W}_{n}\cap\Omega_{\mathbb{S}+n}=\mathcal{W}_{n}$ 
is seen from the fact that
$\psi_{\sigma}(\alpha)>\Omega_{\mathbb{S}+n}$ for $\sigma>\Omega_{\mathbb{S}+n}$.
\\

\noindent
\ref{prp:pi11.0}.\ref{prp:pi11.0.0}.
$\mathcal{W}_{n+1}\cap\Omega_{\mathbb{S}+n+1}\subset W(\mathcal{C}_{n})\subset\mathcal{C}_{n}$ yields $\mathcal{C}_{n+1}\subset\mathcal{C}_{n}$.
\\

\noindent
\ref{prp:pi11.0}.\ref{prp:pi11.0.1}.
Let $\alpha\in\mathcal{W}_{n}$.
Proposition \ref{prp:pi11.0}.\ref{prp:pi11.0.2} yields $\alpha\in\mathcal{W}_{n}=\mathcal{C}_{n}\cap\Omega_{\mathbb{S}+n}$.
Hence $\mathcal{C}_{n}\cap\alpha\subset\mathcal{W}_{n+1} \Rarw \alpha\in\mathcal{W}_{n+1}$.
By induction on $\alpha\in\mathcal{W}_{n}$ we see $\alpha\in\mathcal{W}_{n}\Rarw \alpha\in\mathcal{W}_{n+1}$.
\\

\noindent
\ref{prp:pi11.0}.\ref{prp:pi11.0.3}.
$\mathcal{C}_{N}=W(\mathcal{C}_{N-1})$ is seen from the fact that $\Omega_{\mathbb{S}+N}=\mathbb{K}$
is the largest regular term in $OT_{N}$.
\hspace*{\fill} $\Box$

\begin{proposition}\label{prp:pi11.1}
$\mathcal{C}_{N}\subset G$, where
\\
$\alpha\in G:\Lrarw
\forall n<N(\psi_{\Omega_{\mathbb{S}+n+1}}(\alpha)\in OT_{N} \Rarw \psi_{\Omega_{\mathbb{S}+n+1}}(\alpha)\in\mathcal{W}_{n+1})$.
\end{proposition}
{\it Proof}.\hspace{2mm}
It suffices to show that $G$ is progressive in $\mathcal{C}_{N}$, 
i.e., $Prg[\mathcal{C}_{N},G]$ by 
Proposition \ref{prp:pi11.0}.\ref{prp:pi11.0.3}.
Suppose $\alpha\in\mathcal{C}_{N}$ and $\mathcal{C}_{N}\cap\alpha\subset G$.
We need to show $\alpha\in G$. Let $0<n\leq N$ with $\psi_{\Omega_{\mathbb{S}+n}}(\alpha)\in OT_{N}$.
We show $\psi_{\Omega_{\mathbb{S}+n}}(\alpha)\in\mathcal{W}_{n}$.
By  Proposition \ref{prp:pi11.0}.\ref{prp:pi11.0.0} we obtain $\alpha\in\mathcal{C}_{N}\subset\mathcal{C}_{n-1}$, and hence
$\psi_{\Omega_{\mathbb{S}+n}}(\alpha)\in\mathcal{C}_{n-1}$.
Next we show the following by induction on the length $\ell\gamma$ of ordinal terms $\gamma$.
\[
\forall m\geq n[\gamma\in\mathcal{C}_{m-1}\cap\psi_{\Omega_{\mathbb{S}+m}}(\alpha) \Rarw \gamma\in\mathcal{W}_{m}].
\]
Suppose $\gamma\in\mathcal{C}_{m-1}\cap\psi_{\Omega_{\mathbb{S}+m}}(\alpha)$ for an $m\geq n$.
Propositions \ref{prp:pi11.0}.\ref{prp:pi11.0.2} and \ref{prp:pi11.0}.\ref{prp:pi11.0.1}
yield $\mathcal{C}_{m-1}\cap\Omega_{\mathbb{S}+m-1}=\mathcal{W}_{m-1}\subset\mathcal{W}_{m}$, and
$\Omega_{\mathbb{S}+m-1}\in\mathcal{W}_{m}$.
Hence we can assume $\Omega_{\mathbb{S}+m-1}<\gamma<\psi_{\Omega_{\mathbb{S}+m}}(\alpha)$.
Furthermore by IH and Proposition \ref{prp:pi11.0}.\ref{prp:pi11.0.4} 
we can assume $\gamma=\psi_{\Omega_{\mathbb{S}+m}}(\beta)$ for a $\beta$.
Then $\beta\in\mathcal{C}_{m-1}\cap\alpha$.
By $\mathcal{C}_{N}\cap\alpha\subset G$ it suffices to show $\beta\in\mathcal{C}_{N}$.
Let 
$\Omega_{\mathbb{S}+m-1}<\psi_{\Omega_{\mathbb{S}+k}}(\delta)\in 
k_{\Omega_{\mathbb{S}+m-1}}(\beta)$ for a $k\geq m\geq n$, where
$k_{\Omega_{\mathbb{S}+m-1}}(\beta)$ denotes the set in Definition \ref{df:EGFk}.
We show $\psi_{\Omega_{\mathbb{S}+k}}(\delta)\in\mathcal{W}_{k}$.
We obtain $\delta\in K_{\psi_{\Omega_{\mathbb{S}+k}}(\delta)}(\beta)<\beta<\alpha$, and $\delta\in\mathcal{C}_{k-1}$.
Hence $\psi_{\Omega_{\mathbb{S}+k}}(\delta)\in \mathcal{C}_{k-1}\cap\psi_{\Omega_{\mathbb{S}+k}}(\alpha)$.
IH yields $\psi_{\Omega_{\mathbb{S}+k}}(\delta)\in\mathcal{W}_{k}$.
Hence $\beta\in\mathcal{C}_{N}$.

Thus $\psi_{\Omega_{\mathbb{S}+n}}(\alpha)\in\mathcal{C}_{n-1}$ and
$\mathcal{C}_{n-1}\cap\psi_{\Omega_{\mathbb{S}+n}}(\alpha)\subset\mathcal{W}_{n}$ are shown.
Therefore $\psi_{\Omega_{\mathbb{S}+n}}(\alpha)\in W(\mathcal{C}_{n-1})\cap\Omega_{\mathbb{S}+n}=\mathcal{W}_{n}$.
\hspace*{\fill} $\Box$

\begin{proposition}\label{prp:pi11.2}
$\mathcal{C}^{\mathbb{S}}(\mathcal{W})=\mathcal{C}_{0}\subset\mathcal{W}_{N}$.
\end{proposition}
{\it Proof}.\hspace{2mm}
$\gamma\in\mathcal{C}^{\mathbb{S}}(\mathcal{W}) \Rarw \gamma\in\mathcal{W}_{N}$ is seen by induction on the lengths 
$\ell \gamma$ of ordinal terms $\gamma$ using Proposition \ref{prp:pi11.1}.
\hspace*{\fill} $\Box$

\begin{proposition}\label{prp:psimS}
For $\psi_{\sigma}^{f}(\alpha)\in OT_{N}$,
if $\psi_{\sigma}^{f}(\alpha)\in\mathcal{G}$, then $\psi_{\sigma}^{f}(\alpha)\in \mathcal{W}$.
\end{proposition}
{\it Proof}.\hspace{2mm}
Let $\gamma=\psi_{\sigma}^{f}(\alpha)\in OT_{N}$ and $\gamma\in\mathcal{G}$.
Then $\gamma\in\mathcal{C}^{\gamma}(\mathcal{W})$ and $\mathcal{C}^{\gamma}(\mathcal{W})\cap\gamma\subset\mathcal{W}$.
We obtain $K(f)\subset\mathcal{C}^{\gamma}(\mathcal{W})$.
On the other hand we have $E_{\mathbb{S}}(K(f))\subset\mathcal{C}^{\gamma}(\mathcal{W})\cap\gamma\subset\mathcal{W}$
by (\ref{eq:notationsystem.11}) in Definition \ref{df:notationsystem}.
Hence
$K(f)\subset\mathcal{C}^{\mathbb{S}}(\mathcal{W})$.
We obtain $K(f)\subset\mathcal{W}_{N}$ by Proposition \ref{prp:pi11.2}.

For ordinals $\gamma=\psi_{\sigma}^{f}(\alpha)$ let us associate ordinals $o(\gamma)=o(m(\gamma))=o(f)$
with the ordinal $o(f)<\Gamma_{\mathbb{K}+1}$ in Definition \ref{df:Lam.of}.
From $K(f)\subset\mathcal{W}_{N}$ we see
 $o(\gamma)\in\mathcal{W}_{N}$ for $\gamma=\psi_{\sigma}^{f}(\alpha)\in\mathcal{G}$.
We see $\sigma\prec \rho<\mathbb{S} \Rarw o(\sigma)<o(\rho)$
from Definition \ref{df:notationsystem}.\ref{df:notationsystem.11}
 and Lemma \ref{lem:of}.
Therefore $\gamma\in\mathcal{W}$ is seen by induction on $o(\gamma)\in\mathcal{W}_{N}$ using Lemma \ref{th:3wf16}.
\hspace*{\fill} $\Box$

\begin{lemma}\label{lem:mS}
$\mathbb{S}\in\mathcal{W}$.
\end{lemma}
{\it Proof}.\hspace{2mm}
We have $\mathbb{S}\in\mathcal{C}^{\mathbb{S}}(\mathcal{W})$ and $\mathcal{C}^{\mathbb{S}}(\mathcal{W})\cap\mathbb{S}=\mathcal{C}_{0}\cap\mathbb{S}=\mathcal{W}_{0}\subset\mathcal{W}$.
We obtain $\mathbb{S}\in\mathcal{G}$.
Lemma \ref{th:3wf16} and Proposition \ref{prp:psimS} yield $\mathbb{S}\in\mathcal{W}$.
\hspace*{\fill} $\Box$

\begin{proposition}\label{prp:Wn}
Each $\mathcal{W}_{n}$ is a distinguished set, and $\mathcal{W}_{n}=\mathcal{W}\cap\Omega_{\mathbb{S}+n}$.
\end{proposition}
{\it Proof}.\hspace{2mm}
$\mathcal{W}_{0}=\mathcal{W}\cap\mathbb{S}$ is a distinguished set by Proposition \ref{lem:3wf6}.
Lemma \ref{lem:mS} yields
$\mathcal{W}_{1}=W(\mathcal{C}^{\mathbb{S}}(\mathcal{W}\cap\mathbb{S}))\cap\Omega_{\mathbb{S}+1}=W(\mathcal{C}^{\mathbb{S}}(\mathcal{W}))\cap\Omega_{\mathbb{S}+1}=\mathcal{W}\cap\Omega_{\mathbb{S}+1}$.
Assume that $\mathcal{W}_{n}=\mathcal{W}\cap\Omega_{\mathbb{S}+n}$ for $n>0$.

We have $\Omega_{\mathbb{S}+n}\in\mathcal{G}$, i.e., 
$\mathcal{C}^{\Omega_{\mathbb{S}+n}}(\mathcal{W})\cap\Omega_{\mathbb{S}+n}=\mathcal{C}^{\Omega_{\mathbb{S}+n}}(\mathcal{W}_{n})\cap\Omega_{\mathbb{S}+n}=\mathcal{W}_{n}\subset\mathcal{W}$.
Let $\gamma\in\mathcal{G}$ with $\gamma\prec\Omega_{\mathbb{S}+n}$.
Since there is no $\delta\prec\gamma$, Lemma \ref{th:3wf16} yields that
$\gamma\in W(\mathcal{C}^{\gamma}(\mathcal{W}))$ and $W(\mathcal{C}^{\gamma}(\mathcal{W}))\cap\gamma^{+}$ is a distinguished set,
and hence $\gamma\in\mathcal{W}$.
Again from Lemma \ref{th:3wf16} we see that
$\Omega_{\mathbb{S}+n}\in W(\mathcal{C}^{\Omega_{\mathbb{S}+n}}(\mathcal{W}))\cap\Omega_{\mathbb{S}+n+1}$, and 
$W(\mathcal{C}^{\Omega_{\mathbb{S}+n}}(\mathcal{W}))\cap\Omega_{\mathbb{S}+n+1}$ is a distinguished set.
We conclude that
$\mathcal{W}_{n+1}=W(\mathcal{C}^{\Omega_{\mathbb{S}+n}}(\mathcal{W}_{n}))\cap\Omega_{\mathbb{S}+n+1}=W(\mathcal{C}^{\Omega_{\mathbb{S}+n}}(\mathcal{W}))\cap\Omega_{\mathbb{S}+n+1}=\mathcal{W}\cap\Omega_{\mathbb{S}+n+1}$.
\hspace*{\fill} $\Box$

\begin{proposition}\label{prp:Omega}
Let $\alpha=\Omega_{a}\in OT_{N}$, where either $a=1$ or 
$a=\kappa+n$ with $\kappa\in (\Psi_{N}\cup\{\mathbb{S}\})\cap\mathcal{W}$
and $0<n\leq N$.
Then $\alpha\in\mathcal{W}$.
\end{proposition}
{\it Proof}.\hspace{2mm}
By Propositions \ref{prp:pi11.0}.\ref{prp:pi11.0.2} and \ref{prp:Wn} we can assume either $\alpha=\Omega_{1}$
or $\alpha=\Omega_{\kappa+n}$ with $\kappa\in\Psi_{N}\cap\mathcal{W}$.
Consider the former.
As in Proposition \ref{prp:Wn} we see $\Omega_{1}\in\mathcal{W}$ as follows.

We have $\Omega_{1}\in\mathcal{G}$, i.e., 
$\mathcal{C}^{\Omega_{1}}(\mathcal{W})\cap\Omega_{1}\subset\mathcal{W}$.
Let $\gamma\in\mathcal{G}$ with $\gamma\prec\Omega_{1}$.
Since there is no $\delta\prec\gamma$, Lemma \ref{th:3wf16} yields that
$\gamma\in W(\mathcal{C}^{\gamma}(\mathcal{W}))$ and $W(\mathcal{C}^{\gamma}(\mathcal{W}))\cap\gamma^{+}$ is a distinguished set,
and hence $\gamma\in\mathcal{W}$.
Again from Lemma \ref{th:3wf16} we see that
$\Omega_{1}\in W(\mathcal{C}^{\Omega_{1}}(\mathcal{W}))\cap\Omega_{2}$, and 
$W(\mathcal{C}^{\Omega_{1}}(\mathcal{W}))\cap\Omega_{2}$ is a distinguished set.
Hence $\Omega_{1}\in W(\mathcal{C}^{\Omega_{1}}(\mathcal{W}))\cap\Omega_{2}\subset\mathcal{W}$.
\hspace*{\fill} $\Box$

\begin{definition}\label{df:id4wfA}
{\rm 
For irreducible functions $f$ let
\[
f\in J:\Lrarw K(f)\subset\mathcal{W}.
\]

For $a\in OT_{N}$ and irreducible functions $f$, define:
\begin{eqnarray}
 A(a,f) & :\Lrarw &
 \forall\sigma\in\mathcal{W}[\psi_{\sigma}^{f}(a)\in OT_{N} \Rarw \psi_{\sigma}^{f}(a)\in\mathcal{W}].
\\
\mbox{{\rm MIH}}(a) & :\Lrarw &
 \forall b\in\mathcal{W}\cap a\forall f\in J \, A(b,f).
\\
\mbox{{\rm SIH}}(a,f) & :\Lrarw &
 \forall g\in J [g<^{0}_{lx}f  \Rarw A(a,g)].
\end{eqnarray}
}
\end{definition}

\begin{lemma}\label{th:id5wf21}
Assume $a\in\mathcal{W}$, $f\in J$, $\mbox{{\rm MIH}}(a)$, and $\mbox{{\rm SIH}}(a,f)$ in Definition \ref{df:id4wfA}.
Then
\[
 \forall\kappa\in\mathcal{W}[\psi_{\kappa}^{f}(a)\in OT_{N} \Rarw \psi_{\kappa}^{f}(a)\in\mathcal{W}].
\]
\end{lemma}
{\it Proof}.\hspace{2mm}  
Let $\alpha_{1}=\psi_{\kappa}^{f}(a)\in OT_{N}$ with $a,\kappa\in\mathcal{W}$ and $f\in J$. 
By Proposition \ref{prp:psimS} it suffices to show $\alpha_{1}\in\mathcal{G}$.

By Proposition \ref{lem:3.11.632} we have
$\{\kappa,a\}\cup K(f)\subset\mathcal{C}^{\alpha_{1}}(\mathcal{W})$, and hence $\alpha_{1}\in\mathcal{C}^{\alpha_{1}}(\mathcal{W})$.
It suffices to show the following claim.
\begin{equation}\label{clm:id5wf21.1}
\forall\beta_{1}\in\mathcal{C}^{\alpha_{1}}(\mathcal{W})\cap\alpha_{1}[\beta_1\in\mathcal{W}].
\end{equation}
{\it Proof} of (\ref{clm:id5wf21.1}) by induction on $\ell\beta_1$. 
Assume $\beta_{1}\in\mathcal{C}^{\alpha_{1}}(\mathcal{W})\cap\alpha_{1}$ and let
\[
\mbox{LIH} :\Lrarw
\forall\gamma\in\mathcal{C}^{\alpha_{1}}(\mathcal{W})\cap\alpha_{1}[\ell\gamma<\ell\beta_{1} \Rarw \gamma\in\mathcal{W}].
\]

We show $\beta_1\in\mathcal{W}$.

\noindent
{\bf Case 0}. 
$\beta_1\not\in\mathcal{E}(\beta_{1})$:
Assume $\beta_{1}\not\in\mathcal{W}$.
Then $S(\beta_{1})\subset\mathcal{C}^{\alpha_{1}}(\mathcal{W})\cap\alpha_{1}$. 
LIH yields $S(\beta_{1})\subset\mathcal{W}$. 
Hence we conclude $\beta_{1}\in\mathcal{W}$ from Proposition \ref{lem:3wf8}.
\\

In what follows consider the cases when $\beta_{1}=\psi_{\pi}^{g}(b)$ for some $\pi,b,g$.
We can assume $\{\pi,b\}\cup K(g)\subset\mathcal{C}^{\alpha_{1}}(\mathcal{W})$.
There are several cases according to Propositions \ref{prp:psicomparison} and \ref{prp:l.5.4}.
\\
{\bf Case 1}. $\pi\leq\alpha_{1}$: 
Then $\{\beta_{1}\}=G_{\pi}(\beta_{1})\subset\mathcal{W}$ by $\beta_{1}\in\mathcal{C}^{\alpha_{1}}(\mathcal{W})$ and Proposition \ref{lem:KC}.
\\
{\bf Case 2}.
$b<a$, $\beta_{1}<\kappa$ and $K_{\alpha_{1}}(\{\pi,b\}\cup K(g))<a$:
Let $B$ denote a set of subterms of $\beta_{1}$ defined recursively as follows.
First $\{\pi,b\}\cup K(g)\subset B$.
Let $\alpha_{1}\leq\beta\in B$. 
If $\beta=_{NF}\gamma_{m}+\cdots+\gamma_{0}$, then $\{\gamma_{i}:i\leq m\}\subset B$.
If $\beta=_{NF}\varphi\gamma\delta$, then $\{\gamma,\delta\}\subset B$.
If $\beta=_{NF}\Omega_{\gamma}$, then $\gamma\in B$.
If $\beta=\psi_{\sigma}^{h}(c)$, then $\{\sigma,c\}\cup K(h)\subset B$.

Then from $\{\pi,b\}\cup K(g)\subset\mathcal{C}^{\alpha_{1}}(\mathcal{W})$ we see inductively that
$B\subset\mathcal{C}^{\alpha_{1}}(\mathcal{W})$.
Hence by LIH we obtain $B\cap\alpha_{1}\subset\mathcal{W}$.
Moreover if $\alpha_{1}\leq\psi_{\sigma}^{h}(c)\in B$, then $c\in K_{\alpha_{1}}(\{\pi,b\}\cup K(g))<a$.
We claim that
\begin{equation}\label{eq:case2A}
\forall\beta\in B(\beta\in\mathcal{W})
\end{equation}
{\it Proof} of (\ref{eq:case2A}) by induction on $\ell\beta$.
Let $\beta\in B$. We can assume that $\alpha_{1}\leq\beta=\psi_{\sigma}^{h}(c)$ by induction hypothesis on the lengths.
Then by induction hypothesis we have $\{\sigma,c\}\cup K(h)\subset\mathcal{W}$.
On the other hand we have $c<a$.
$\mbox{MIH}(a)$ yields $\beta\in\mathcal{W}$.
Thus (\ref{eq:case2A}) is shown.
\hspace*{\fill} $\Box$
\\

In particular we obtain $\{\pi,b\}\cup K(g)\subset \mathcal{W}$.
Moreover we have $b<a$.
Therefore once again $\mbox{MIH}(a)$ yields $\beta_{1}\in\mathcal{W}$.
\\
{\bf Case 3}. 
$b=a$, $\pi=\kappa$, $\forall\delta\in K(g)(K_{\alpha_{1}}(\delta)<a)$ and $g<^{0}_{lx}f$: 
As in (\ref{eq:case2A}) we see that $K(g)\subset\mathcal{W}$ from $\mbox{MIH}(a)$.
$\mbox{SIH}(a,f)$ yields $\beta_{1}\in\mathcal{W}$.
\\
{\bf Case 4}.
$a\leq b\leq K_{\beta_{1}}(\delta)$ for some $\delta\in K(f)\cup\{\kappa,a\}$:
It suffices to find a $\gamma$ such that $\beta_{1}\leq\gamma\in\mathcal{W}\cap\alpha_{1}$.
Then $\beta_{1}\in\mathcal{W}$ follows from $\beta_{1}\in\mathcal{C}^{\alpha_{1}}(\mathcal{W})$ and Proposition \ref{prp:updis}.

$k_{\delta}(\alpha)$ denotes the set in Definition \ref{df:EGFk}.
In general we see that $a\in K_{\delta}(\alpha)$ iff $\psi_{\sigma}^{h}(a)\in k_{\delta}(\alpha)$ for some $\sigma,h$,
and for each $\psi_{\sigma}^{h}(a)\in k_{\delta}(\psi_{\sigma_{0}}^{h_{0}}(a_{0}))$ there exists a sequence
$\{\alpha_{i}\}_{i\leq m}$ of subterms of $\alpha_{0}=\psi_{\sigma_{0}}^{h_{0}}(a_{0})$ such that 
$\alpha_{m}=\psi_{\sigma}^{h}(a)$, 
$\alpha_{i}=\psi_{\sigma_{i}}^{h_{i}}(a_{i})$ for some $\sigma_{i},a_{i},h_{i}$,
and for each $i<m$,
$\delta\leq\alpha_{i+1}\in\mathcal{E}(C_{i})$ for $C_{i}=\{\sigma_{i},a_{i}\}\cup K(h_{i})$.

Let $\delta\in K(f)\cup\{\kappa,a\}$ such that $b\leq \gamma$
for a $\gamma\in K_{\beta_{1}}(\delta)$.
Pick an $\alpha_{2}=\psi_{\sigma_{2}}^{h_{2}}(a_{2})\in \mathcal{E}(\delta)$ 
such that $\gamma\in K_{\beta_{1}}(\alpha_{2})$, and 
an $\alpha_{m}=\psi_{\sigma_{m}}^{h_{m}}(a_{m})\in k_{\beta_{1}}(\alpha_{2})$ for some $\sigma_{m},h_{m}$ and
$a_{m}\geq b\geq a$.
We have $\alpha_{2}\in\mathcal{W}$ by $\delta\in\mathcal{W}$.
If $\alpha_{2}<\alpha_{1}$, then
$\beta_{1}\leq\alpha_{2}\in\mathcal{W}\cap\alpha_{1}$, and we are done.
Assume $\alpha_{2}\geq\alpha_{1}$.
Then $a_{2}\in K_{\alpha_{1}}(\alpha_{2})<a\leq b$, and $m>2$.

Let $\{\alpha_{i}\}_{2\leq i\leq m}$ be the sequence of subterms of $\alpha_{2}$ such that
$\alpha_{i}=\psi_{\sigma_{i}}^{h_{i}}(a_{i})$ for some $\sigma_{i},a_{i},h_{i}$,
and for each $i<m$,
$\beta_{1}\leq\alpha_{i+1}\in\mathcal{E}(C_{i})$ for $C_{i}=\{\sigma_{i},a_{i}\}\cup K(h_{i})$.

Let $\{n_{j}\}_{0\leq j\leq k}\, (0<k\leq m-2)$ be the increasing sequence $n_{0}<n_{1}<\cdots<n_{k}\leq m$ 
defined recursively by $n_{0}=2$, and assuming $n_{j}$ has been defined so that
$n_{j}<m$ and $\alpha_{n_{j}}\geq\alpha_{1}$, $n_{j+1}$ is defined by
$n_{j+1}=\min(\{i: n_{j}< i<m, \alpha_{i}<\alpha_{n_{j}}\}\cup\{m\})$.
If either $n_{j}=m$ or $\alpha_{n_{j}}<\alpha_{1}$, then $k=j$ and $n_{j+1}$ is undefined.
Then we claim that
\begin{equation}\label{eq:case4A}
\forall j\leq k(\alpha_{n_{j}}\in\mathcal{W}) \spand \alpha_{n_{k}}<\alpha_{1}
\end{equation}
{\it Proof} of (\ref{eq:case4A}).
By induction on $j\leq k$ we show first that $\forall j\leq k(\alpha_{n_{j}}\in\mathcal{W})$. 
We have $\alpha_{n_{0}}=\alpha_{2}\in\mathcal{W}$.
Assume $\alpha_{n_{j}}\in\mathcal{W}$ and $j<k$.
Then $n_{j}<m$, i.e., $\alpha_{n_{j+1}}<\alpha_{n_{j}}$, and 
by $\alpha_{n_{j}}\in\mathcal{C}^{\alpha_{n_{j}}}(\mathcal{W})$, we have $C_{n_{j}}\subset\mathcal{C}^{\alpha_{n_{j}}}(\mathcal{W})$,
and hence $\alpha_{n_{j}+1}\in\mathcal{E}(C_{n_{j}})\subset\mathcal{C}^{\alpha_{n_{j}}}(\mathcal{W})$.
We see inductively that
$\alpha_{i}\in \mathcal{C}^{\alpha_{n_{j}}}(\mathcal{W})$ for any $i$ with $n_{j}\leq i\leq n_{j+1}$.
Therefore $\alpha_{n_{j+1}}\in \mathcal{C}^{\alpha_{n_{j}}}(\mathcal{W})\cap\alpha_{n_{j}}\subset\mathcal{W}$ by Proposition \ref{prp:updis}.

Next we show that $\alpha_{n_{k}}<\alpha_{1}$.
We can assume that $n_{k}=m$.
This means that $\forall i(n_{k-1}\leq i<m \Rarw \alpha_{i}\geq\alpha_{n_{k-1}})$.
We have
$\alpha_{2}=\alpha_{n_{0}}>\alpha_{n_{1}}>\cdots>\alpha_{n_{k-1}}\geq\alpha_{1}$, and
$\forall i<m(\alpha_{i}\geq\alpha_{1})$.
Therefore $\alpha_{m}\in k_{\alpha_{1}}(\alpha_{2})\subset k_{\alpha_{1}}(\{\kappa,a\}\cup K(h))$, i.e.,
$a_{m}\in K_{\alpha_{1}}(\{\kappa,a\}\cup K(h))$ for $\alpha_{m}=\psi_{\sigma_{m}}^{h_{m}}(a_{m})$.
On the other hand we have $K_{\alpha_{1}}(\{\kappa,a\}\cup K(h))<a$ for 
$\alpha_{1}=\psi_{\sigma}^{h}(a)$.
Thus $a\leq a_{m}<a$, a contradiction.

(\ref{eq:case4A}) is shown, and we obtain $\beta_{1}\leq\alpha_{n_{k}}\in\mathcal{W}\cap\alpha_{1}$.

This completes a proof of (\ref{clm:id5wf21.1}) and of the lemma.
\hspace*{\fill} $\Box$

\begin{lemma}\label{lem:psiw}
For $\psi_{\kappa}^{f}(a)\in OT_{N}$, if
$\{\kappa,a\}\cup K(f)\subset\mathcal{W}$, then $\psi_{\kappa}^{f}(a)\in \mathcal{W}$.
\end{lemma}
{\it Proof}.\hspace{2mm}
$\forall a\in\mathcal{W}\forall f\in J\forall \kappa\in\mathcal{W}[
\psi_{\kappa}^{f}(a)\in OT_{N} \Rightarrow \psi_{\kappa}^{f}(a)\in\mathcal{W}]$
is seen by main induction on $a\in\mathcal{W}$ with subsidiary induction
on $f\in J$ using Lemma \ref{th:id5wf21} and Lemma \ref{lem:oflx}.
Note that the ordinal $o(f)\in\mathcal{W}$ in the proof of Proposition \ref{prp:psimS}
 if $f\in J$, i.e., if $K(f)\subset\mathcal{W}$.
\hspace*{\fill} $\Box$
\\

\noindent
{\it Proof} of Theorem \ref{th:wf}.
Recall that $\mathcal{W}$ is a distinguished set by Proposition \ref{lem:3wf6}, and
distinguished sets are well founded.
By Proposition \ref{lem:3wf8}, $\mathcal{W}$ is closed under $+,\varphi$ and $0\in\mathcal{W}$.
By Proposition \ref{prp:Omega}, $\{\Omega_{1}\}\cup\{\Omega_{\mathbb{S}+n}: 0<n\leq N\}\subset\mathcal{W}$,
and $\kappa\in\Psi_{N}\cap\mathcal{W}\Rarw \Omega_{\kappa+n}\in\mathcal{W}$.
By Lemma \ref{lem:mS}, $\mathbb{S}\in\mathcal{W}$, and $\mathcal{W}$ is closed under 
$(\kappa,a,f)\mapsto \psi_{\kappa}^{f}(a)$ by Lemma \ref{lem:psiw}.
Hence we see that $\alpha\in OT_{N}\Rarw\alpha\in\mathcal{W}$ by induction on $\ell\alpha$.

\end{document}